\documentclass[final,12pt]{elsarticle}

\usepackage{amssymb}
\usepackage{amsmath}
\usepackage{amsthm}

\usepackage{comment}

\usepackage{booktabs}

\usepackage[dvipsnames]{xcolor}

\usepackage{algorithm}
\usepackage{algpseudocode}

\usepackage{subfig}
\usepackage{caption}

\usepackage{cases}

\usepackage{stmaryrd}

\newcommand{\vertiii}[1]{{\left\vert\kern-0.25ex\left\vert\kern-0.25ex\left\vert #1 
    \right\vert\kern-0.25ex\right\vert\kern-0.25ex\right\vert}}

\newcommand{\dm}[1]{\textcolor{blue}{#1}}

\theoremstyle{definition}
\newtheorem{lemma}{Lemma}[section]
\newtheorem{remark}{Remark}[section]
\newtheorem{proposition}{Proposition}[section]
\journal{Computers \& Mathematics with Applications}

\begin{document}

\begin{frontmatter}



\title{A novel Mortar Method Integration using Radial Basis Functions}


\affiliation[UniPD]{organization={Dept. of Civil, Environmental and Architectural Engineering, University of Padova},
            state={Padova},
            country={Italy}}

\author[UniPD]{Daniele Moretto\corref{coraut}} 
\ead{daniele.moretto.3@phd.unipd.it}
\author[UniPD]{Andrea Franceschini}
\author[UniPD]{Massimiliano Ferronato}

\cortext[coraut]{Corresponding author}

\begin{abstract}
The growing availability of computational resources has significantly increased the interest of the scientific community in performing complex multi-physics and multi-domain simulations. However, the generation of appropriate computational grids for such problems often remains one of the main bottlenecks. The use of a domain partitioning with non-conforming grids is a possible solution, which, however, requires the development of robust and efficient inter-grid interpolation operators to transfer a scalar or a vector field from one domain to another. This work presents a novel approach for interpolating quantities across non-conforming meshes within the framework of the classical mortar method, where weak continuity conditions are enforced. The key contribution is the introduction of a novel strategy that uses mesh-free Radial Basis Function (RBF) interpolations to compute the mortar integral, offering a compelling alternative to traditional projection-based methods. We propose an efficient algorithm tailored for complex three-dimensional settings allowing for potentially significant savings in the overall computational cost and ease of implementation, with no detrimental effects on the numerical accuracy. The formulation, analysis, and validation of the proposed RBF-based algorithm is discussed with the aid of a set of numerical examples, demonstrating its effectiveness. Furthermore, the details of the implementation are discussed and a test case involving a complex geometry is presented, to illustrate the applicability and advantages of our approach in real-world problems.
\end{abstract}



\begin{keyword}
Mortar Method \sep Inter-grid Interpolation \sep Radial Basis Functions (RBF) \sep Multi-physics and Multi-domain Simulations




\end{keyword}
\end{frontmatter}



\section{Introduction}

Over the last decades, the growing availability of computing capabilities has increased the interest in the numerical simulation of complex multi-physics and multi-domain problems. 
For instance, similar frameworks arise when combining multi-scale processes in subsurface modeling \cite{castelletto2017multiscale,mehmani2019multiscale,cusini2021simulation,cusini2022field,franceschini2020algebraically}. 
In this context, splitting the overall domain into different parts with independent spatial discretizations offers several well-known attractive advantages \cite{belgacem1999mortar}. Using also different temporal discretizations  \cite{brun2015two,gravouil2015heterogeneous} is another appealing option that will not be considered in this paper. Focusing on the spatial approach, it requires the definition of an inter-grid interpolation operator, which transfers a certain variable field across different discretizations. 
Such an operator acts on the lower-dimensional interface separating a pair of adjacent domains, which are generally characterized by two non-conforming grids. The components of this pair of lower-dimensional objects in contact each other are typically referred to as \textit{master} and \textit{slave} sides, respectively.

Over the years many numerical strategies have been proposed to compute the inter-grid interpolation operator and an exhaustive review is beyond the scope of this work. Following \cite{de2007review,coniglio2019weighted}, most of these strategies can be grouped into two categories. The first one includes the \textit{collocation methods}, known also as \textit{node-to-segment} approaches within the contact mechanics community, which establish a linear relationship between the degrees of freedom of the slave side and the master side in a strong way. Such relationships may derive from simple point-wise matching conditions or from more accurate projection/interpolation techniques, and must guarantee the global conservation of physical quantities across the interface. Collocation methods are sub-optimal \cite{bernardi1993domain} and their accuracy strongly depends on the choice of the master and the slave sides. The second category, denoted \textit{weighted residual methods} \cite{de2007review} and also known as \textit{segment-to-segment} approaches \cite{puso2004mortar,de2021displacement}, includes the techniques where the continuity of quantities at the interface is enforced in a weak form. The \textit{mortar method} \cite{bernardi1993domain,belgacem1997mortar} is perhaps the most popular representative of this family, where continuity is prescribed by introducing a set of Lagrange multipliers on the slave interface \cite{belgacem1999mortar}.

\begin{figure}
    \centering
    \includegraphics[width=0.6\linewidth]{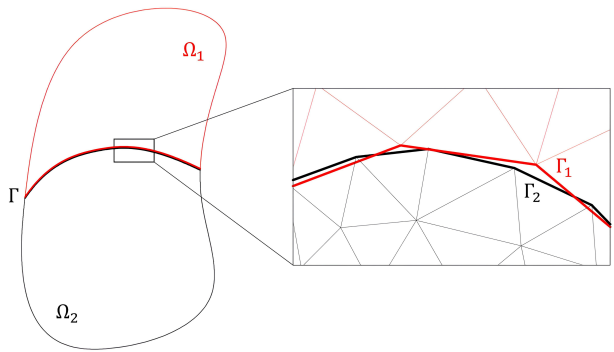}
    \caption{Geometrically non-conforming grids arising from curved interface $\Gamma$, $d=2$.}
    \label{fig_geomNonConf}
\end{figure}

In this work, we focus on the mortar method applied to a general second-order elliptic differential problem discretized by finite elements. For the sake of simplicity and with no loss of generality for the scopes of the present discussion, we resort to the following preliminary assumptions: 
\begin{itemize}
    \item homogeneous global boundary conditions;
    \item a common interface $\Gamma$ is shared by two domains only. 
\end{itemize}
With the second hypothesis, we exclude the presence of cross points appearing in more than two interfaces, whose treatment involves some technicalities \cite{puso2003mesh}.
Conversely, no assumptions are made on the geometry of the interface. Curved interfaces may generally lead to situations like the one shown in Figure \ref{fig_geomNonConf}, where the grids are said to be \textit{geometrically non-conforming}.

Let $\Omega \subset \mathbb{R}^d$, $d = 2, \,3$, be an open domain and $\partial \Omega$ its Lipschitz boundary, decomposed into Dirichlet $\partial \Omega_D$ and Neumann $\partial \Omega_N$ subsets. 
Define on $\Omega$ a partition into a master subset $\Omega_1$ and a slave subset $\Omega_2$, such that $\Omega \equiv \Omega_1 \cup \Omega_2$ and $\Gamma \equiv \partial \Omega_1 \cap \partial\Omega_2$ is their common interface. Let $u$ be an appropriate continuous function defined onto $\Omega$, with
$u_1$ and $u_2$ the restrictions of $u$ to each subdomain $\Omega_1$ and $\Omega_2$, respectively.
Given $f \in L^2(\Omega)$, we aim at finding the function $u$ that solves the \textit{transmission problem}
\begin{align}
    \begin{cases}
        L_1 u_1 = f  &\text{in $\Omega_1$}, \\
        L_2 u_2 = f  &\text{in $\Omega_2$}, \\
        u_1 = u_2 &\text{on $\Gamma$ },\\
        \partial_{L_1} u_1 + \partial_{L_2} u_2 = 0 \ &\text{on $\Gamma$ }.
    \end{cases}
    \label{Transmission problem}
\end{align}
In equation \eqref{Transmission problem}, $L_k$ is a generic second-order elliptic differential operator and $\partial_{L_k}$ its conormal derivative. We introduce the following spaces:
\begin{align}
    \mathcal{X}_k &= \{ v_k \in H^1(\Omega_k) \ | \ v_k = 0 \text{ on } \partial\Omega_D \cap \partial \Omega_k\}, \quad k=1,2, \nonumber\\
    \mathcal{X} &= \{ v \in L^2(\Omega)  \ | \ v \vert_{\Omega_k} = v_k \in \mathcal{X}_k  \ \forall k\},  \label{Continuous spaces}\\ 
    \mathcal{M} &= \{ \mu \in H^{-1/2}(\Gamma)\}. \nonumber    
\end{align}
$\mathcal M$ is a space of Lagrange multipliers representing the flux of the solution across the interface $\Gamma$. Let also $a_k(u_k,v_k)$ be the local bilinear form associated to $L_k$, with the global one defined as $a(u,v) = \sum_{k=1,2} a_k(u_k,v_k)$.
Then, the mortar formulation with Lagrange multipliers associated to the transmission problem \eqref{Transmission problem} reads:
Find $\{u, \, \lambda\} \ \in \ \{ \mathcal{X} \times \mathcal{M} \}$ such that
\begin{subequations}
\begin{align}
    a(u,v) + b(\lambda,v) &= f(v) &\forall v \in \mathcal{X}, \\ 
            b(u,\mu) &= 0 &\forall \mu \in \mathcal{M},
\label{Mortar Lagrange multipliers}
\end{align}
\label{Continuous mortar form}
\end{subequations}
with 
\begin{align}
    b(u,\mu) = \langle \mu, \llbracket u \rrbracket \rangle_{\Gamma}.
\label{Mortar integral}
\end{align}
In equation \eqref{Mortar integral}, $\llbracket u \rrbracket = (u_2 - u_1)\vert_\Gamma$ denotes the jump of $u$ across the interface $\Gamma$.
Moving to the discrete setting, we consider independent partitions $\mathcal T_k$ of $\Omega_k$, consisting of either triangular/quadrilateral elements in 2D or tetrahedral/hexahedral elements in 3D with characteristic mesh size $h_k$.
We denote by $\Gamma_k$ 
the discretization of $\Gamma$ 
induced by each partition $\mathcal{T}_k$ and by $\mathcal{E}_k$ the collection of edges/faces constituting $\Gamma_k$.
The discrete counterparts of the spaces \eqref{Continuous spaces} read
\begin{align}
    \mathcal{X}_{h,k} &= \{ v_{h,k} \in \mathcal{X}_k \ | \ v_{h,k}\vert_T \in \mathbb{P}_p, \ \forall T \in \mathcal{T}_{k} \}, \quad k=1,2, \nonumber\\
    \mathcal{X}_h &= \{ v_h \in \mathcal{X}, \ | \  \ v_h \vert_{\Omega_k} = v_{h,k} \in \mathcal{X}_{h,k} \ \forall k \}, \\
    \mathcal{M}_h &= \{ \mu_h \in L^2(\Gamma_2), \ | \ \mu_h \vert_E   \in \mathbb{P}_p, \ \forall E \in \mathcal{E}_{2}\}, \nonumber
    \label{Discrete spaces}
\end{align}
where $\mathbb{P}_p$ denotes the space of polynomials of total degree $\leq p$. In this work, we will focus on $p \leq 2$. 
If we set $u_h = [u_{h,1}, \, u_{h,2}]$, the discrete finite element formulation of the mortar problem \eqref{Continuous mortar form} reads: Find $\{u_h, \, \lambda_h \} \in \{ \mathcal{X}_h \times \mathcal{M}_h \}$ such that:
\begin{subequations}
\begin{align}
    a_h(u_h,v_h) + b_h(\lambda_h,v_h) &= f(v_h) &\forall v_h \in \mathcal{X}_h, \\ 
            b_h(u_h,\mu_h) &= 0 &\forall \mu_h \in \mathcal{M}_h, 
\label{Discrete mortar constraint}
\end{align}    
\label{discrete mortar form}
\end{subequations}
where the interface constraint at the discrete level is
\begin{align}
    b_h(u_h,\mu_h) = \int_{\Gamma_2} \mu_h (u_{h,2} - \Pi \, u_{h,1}) \, d\Gamma.
    \label{mortar integral}
\end{align}
In equation \eqref{mortar integral}, a suitable mapping $\Pi: L^2(\Gamma_1) \rightarrow L^2(\Gamma_2)$ must be introduced to deal with the general case of non-matching interfaces, usually defined as a geometric projection of $u_{h,1}$ onto the slave side.
Integral \eqref{mortar integral} is commonly referred to as \textit{mortar integral}. Its computation is the main practical challenge of the mortar method, since the integrand function contains contributions living on different grids.
Though exact integration of \eqref{mortar integral} yields optimal convergence, its implementation in a general 3D setting can be very expensive, since it involves the projection of 3D discretized surfaces one over the other. Hence, in practical applications an inexact integration is typically used to compute \eqref{mortar integral}, even though the mortar optimality might be lost \cite{falletta2007approximate}. 

In this work, we introduce a novel algorithm to compute the mortar integral. The key idea relies on replacing the geometric projection $\Pi$ with a mesh-free interpolation based on Radial Basis Functions (RBF), see \cite{buhmann2000radial,wendland2004scattered,fasshauer2007meshfree} for an overview on this topic. By combining the inexact quadrature of \eqref{mortar integral} with RBF, we obtain an efficient integration algorithm easy to implement also in complex 3D settings.  
Various inter-grid methods based on RBF have been already proposed in \cite{cordero2014radial,deparis2014rescaled,biancolini2018balanced}, just to cite a few. 
The INTERNODES method \cite{deparis2016internodes,gervasio2018analysis} is another optimal non-conforming approach where cross-grid integration is avoided. Here, Rescaled Localized Radial Basis Functions (RL-RBF) \cite{deparis2014rescaled} were suggested to interpolate a variable in the presence of geometrically non-conforming interfaces. However, as compared to the mortar method, INTERNODES does not preserve symmetry,
and total forces and work at the interface are not exactly zero, but vanish optimally with the polynomial degree of the adopted discretization \cite{deparis2022conservation}. 
In contrast to the approaches mentioned above, the proposed algorithm uses RBF directly within the mortar framework. 

The paper is organized as follows: in Section \ref{sect2:DiscreteMortar}, we derive the algebraic form of the mortar method and revise the state-of-the-art quadrature schemes to evaluate the mortar integral. In Section \ref{sect3:RBFalgorithm}, the novel RBF-based algorithm is presented, along with some implementation guidelines. The method is first validated in Section \ref{sect4:NumericalResults}, then tested in a multi-physics and multi-domain application in Section \ref{sect5:Application}. Finally, some concluding remarks close the paper.

\section{Discrete form of the mortar method} \label{sect2:DiscreteMortar}

Let $\mathcal{N}_{\Gamma_k}$ be the set of indices of the nodes of $\Omega_k$ belonging to $\Gamma_k$, and $\mathcal{N}_k$ the set of indices of the remaining nodes of $\Omega_k$, with $\Omega_1$ the master subset and $\Omega_2$ the slave subset. 
By $\{N^k_i\}_{i \in \mathcal{N}_k \cup \mathcal{N}_{\Gamma_k}}$ we denote the standard nodal basis functions for $\mathcal X_{h,k}$, with $\{N^{\Gamma_k}_j\}_{j \in \mathcal{N}_{\Gamma_k} }$ the respective traces over $\mathcal{E}_k$. The traces are 1D shape functions defined over edges and 2D shape functions over the faces for $d=2$ and $d=3$, respectively.
Likewise, $\{\psi_i\}_{i \in \mathcal{N}_{\Gamma_2}}$ denote the basis functions for Lagrange multipliers on the slave side.
The discrete approximations for the introduced variable fields can be expressed as
\begin{align}
    u_{h,k}(\boldsymbol x) = \sum_{i \in \mathcal{N}_k}  N^k_i(\boldsymbol x) \, u_{i,k} &+  \sum_{i \in \mathcal{N}_{\Gamma_k}} N^k_i(\boldsymbol x) \, u_{i,\Gamma_k} , \nonumber \\  
    \lambda_{h}(\boldsymbol x) = \sum_{i \in \mathcal{N}_{\Gamma_2}} \psi_i(\boldsymbol x) \, \lambda_{i} \ \ \ \ \ \ \ \ \ \ 
     &u_{h,k}(\boldsymbol x)\vert_{\Gamma_k} = \sum_{i \in \mathcal{N}_{\Gamma_k}}  N^{\Gamma_k}_i(\boldsymbol x) \, u_{i,\Gamma_k} \label{FEM approximations}
\end{align}
The nodal unknowns $\{ u_{i,k} \}, \, \{ u_{i,\Gamma_k} \}$ and $\{ \lambda_{i} \}$ in equation \eqref{FEM approximations} are collected in the vectors $\boldsymbol{u}_k$, $\boldsymbol{u}_{\Gamma_k}$ and $\boldsymbol{\lambda}$, respectively. Then, the matrix form of the discrete variational problem \eqref{discrete mortar form} reads
\begin{align}
\begin{bmatrix}
    A_{11} & 0 & A_{1\Gamma_1} & 0 & 0 \\
    0 & A_{22} & 0 & A_{2\Gamma_2} & 0 \\
    A_{1\Gamma_1}^T & 0 & A_{\Gamma_1\Gamma_1} & 0 & -S^T \\
    0 & A_{2\Gamma_2}^T & 0 & A_{\Gamma_2\Gamma_2} & D^T \\
    0 & 0 & -S & D & 0 \\
\end{bmatrix}
\begin{bmatrix}
    \boldsymbol{u}_1 \\ \boldsymbol{u}_{2} \\ \boldsymbol{u}_{\Gamma_1} \\ \boldsymbol{u}_{\Gamma_2} \\ \boldsymbol{\lambda}
\end{bmatrix} = \begin{bmatrix}
    \boldsymbol{f}_1 \\ \boldsymbol{f}_{2} \\ \boldsymbol{f}_{\Gamma_1} \\ \boldsymbol{f}_{\Gamma_2} \\ \boldsymbol 0 
\end{bmatrix}
\label{saddlePoint_system}
\end{align}
Matrices $A_{(\cdot,\cdot)}$ and forcing vectors $\boldsymbol f_{(\cdot)}$ are built in the standard way from the respective forms $a_h$ and $f$, whereas matrices $S$ and $D$ represent the coupling of master and slave traces:
\begin{align}
    D[i,j] =  \int_{\Gamma_2}  \psi_i N^{\Gamma_2}_j \; d\Gamma,  \ \ \ \ \ &i,j = 1,\dots,\vert\mathcal{N}_{\Gamma_2}\vert, \label{def_D}\\
    S[i,l] =  \int_{\Gamma_2}  \psi_i \,  \Pi N^{\Gamma_1}_l \; d\Gamma, \ \ \ \ \ &i = 1,\dots,\vert\mathcal{N}_{\Gamma_2}\vert, \, l = 1,\dots, \vert \mathcal{N}_{\Gamma_1} \vert \label{def_M}.
\end{align}
If the Lagrange multiplier basis functions $\psi_i(x)$ are standard polynomials, then $D$ is a classical mass matrix evaluated at the slave interface.  Other choices, as well as general criteria to construct an appropriate set of multipliers, can be found in the literature \cite{belgacem1999mortar,wohlmuth2001iterative}, with the most notable alternative represented by the family of \textit{dual Lagrange multipliers} \cite{wohlmuth2000mortar,puso20043d,flemisch2007stable,popp2010dual,popp2012dual}. 
%
By introducing the interpolation operator $E$ between the master and slave side 
\begin{align}
    E = D^{-1} \, S,
\end{align}
the last equation of system \eqref{saddlePoint_system} can be written as
\begin{align}
    \boldsymbol{u}_{\Gamma_2} = E \, \boldsymbol {u}_{\Gamma_1},
    \label{discrete_mortar_constr_alg}
\end{align}
allowing for a static condensation process. 
Eliminating
both $\boldsymbol u_{\Gamma_2}$ and $\boldsymbol \lambda$ 
yields  
\begin{align}
    \begin{bmatrix}
    A_{11} & 0 & A_{1\Gamma_1} \\
    0 & A_{22} &  A_{2\Gamma_2}E \\
    A_{1\Gamma_1}^T & E^{T}A_{2\Gamma_2}^T & A_{\Gamma_1\Gamma_1} +  E^{T} A_{\Gamma_2\Gamma_2} E 
\end{bmatrix}
\begin{bmatrix}
    \boldsymbol{u}_1 \\ \boldsymbol{u}_{2} \\ \boldsymbol{u}_{\Gamma_1}
\end{bmatrix} = \begin{bmatrix}
    \boldsymbol f_1 \\ \boldsymbol f_{2} \\ \boldsymbol f_{\Gamma_1} + E^{T} \boldsymbol f_{\Gamma_2}
\end{bmatrix}
\label{condensated_system}
\end{align}
which 
preserves symmetry and definiteness of the global problem, if any. 
When adopting \textit{dual Lagrange multipliers}, matrix $D$ is diagonal and the static condensation procedure becomes trivial. 

\subsection{Consistency of the discrete mortar constraint} \label{sect:consistency_mortar}
A fundamental requirement for the consistency of the discrete mortar constraint \eqref{Discrete mortar constraint} 
concerns the properties that the projection of the master basis functions onto the slave side must satisfy, as stated by the result that follows \cite{puso20043d}.

\begin{proposition}
\label{th:paruni}
\emph{
If the discrete mortar constraint \eqref{Discrete mortar constraint} is consistent, then the projection $\Pi$ must be such that:}
\begin{equation}
    \sum_{l=1}^{\vert \mathcal{N}_{\Gamma_1} \vert} \Pi N_l^{\Gamma_1} = \sum_{j=1}^{\vert \mathcal{N}_{\Gamma_2} \vert} N_j^{\Gamma_2}.
    \label{eq:prop_paruni}
\end{equation}
\end{proposition}

\begin{proof}
The consistency of operator $E$ in \eqref{discrete_mortar_constr_alg} implies that constant functions must be preserved in their transfer from the master to the slave interface. Hence, if all the components of $\boldsymbol{u}_{\Gamma_1}$ and $\boldsymbol{u}_{\Gamma_2}$ are equal to any non-zero constant $c\in\mathbb{R}$, it follows that:
\begin{align}
     \sum_{l = 1}^{\vert \mathcal{N}_{\Gamma_1} \vert} E[i,l] = 1 \ \ \ \ \ \ \   \forall \; i=1,\dots,\vert \mathcal{N}_{\Gamma_2} \vert,
\end{align}
or equivalently
\begin{align}
    \sum_{j = 1}^{\vert \mathcal{N}_{\Gamma_2} \vert} D[i,j] = \sum_{l = 1}^{\vert \mathcal{N}_{\Gamma_1} \vert} S[i,l] \ \ \ \ \ \ \ \forall \; i=1,\dots,\vert \mathcal{N}_{\Gamma_2} \vert.
\end{align}
Replacing $D$ and $S$ with their definitions \eqref{def_D} and \eqref{def_M} yields:
\begin{align*}
    \sum_{j = 1}^{\vert \mathcal{N}_{\Gamma_2} \vert}\int_{\Gamma_2}  \psi_i N^{\Gamma_2}_j \; d\Gamma =  \sum_{l = 1}^{\vert \mathcal{N}_{\Gamma_1} \vert} \int_{\Gamma_2}  \psi_i \,  \Pi N^{\Gamma_1}_l \; d\Gamma \, , \ \ \ \ \  \forall \; i= 1,\dots,\vert\mathcal{N}_{\Gamma_2}\vert,
\end{align*}
that can be rearranged as
 \begin{align}
    \int_{\Gamma_2}  \psi_i \sum_{j = 1}^{\vert \mathcal{N}_{\Gamma_2} \vert} N^{\Gamma_2}_j \; d\Gamma =   \int_{\Gamma_2}  \psi_i \sum_{k = 1}^{\vert \mathcal{N}_{\Gamma_1} \vert}\,  \Pi N^{\Gamma_1}_k \; d\Gamma \, ,  \ \ \ \ \  \forall \; i= 1,\dots,\vert\mathcal{N}_{\Gamma_2}\vert.
    \label{proof_PU}
\end{align}
Equation \eqref{proof_PU} holds if the condition \eqref{eq:prop_paruni} is satisfied.        
\end{proof}

\begin{remark}
\label{rem:unity}
The requirement of Proposition \ref{th:paruni} is certainly satisfied if the slave basis functions and the projections of the master basis functions onto the slave side are a partition of unity. For instance, this condition is trivially fulfilled by employing Lagrangian basis functions, as well as by the usual dual basis functions \cite{popp2012mortar}.
\end{remark}

\begin{remark}
The integrals in \eqref{proof_PU} defining the entries of matrices $D$ and $S$ are usually computed numerically. If the slave basis functions and the projections of the master basis functions are partition of unity, then it is necessary to use the same quadrature rule for computing the entries of $D$ and $S$ \cite{puso2004mortar, farah2018mortar}. If not, the necessary condition \eqref{proof_PU} implied by the operator consistency is no longer numerically met.   
\end{remark}

\subsection{Mortar integral computation}
\begin{figure}
    \centering
    \includegraphics[width=0.7\linewidth]{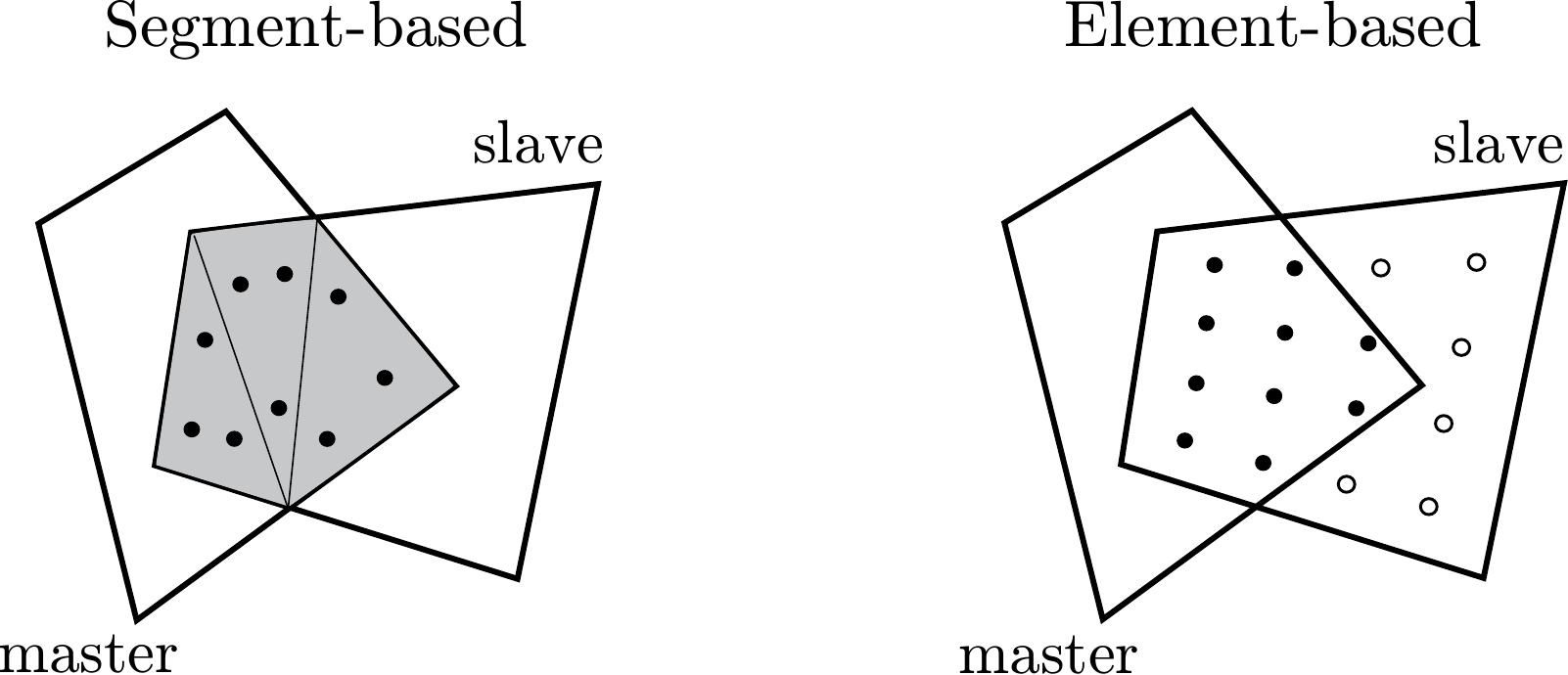}
    \caption{Sketch of the commonly employed strategies to evaluate mortar integral: segment-based (left) and element-based (right).}
    \label{fig:integration_scheme}
\end{figure}

The computation of $D$ does not pose any particular difficulty.
The actual challenge for the implementation of the mortar method relies on the computation of $S$ (equation \eqref{def_M}), since it requires the integration of products of functions defined on different meshes.
Two different approaches are usually employed to compute $S$, namely the \textit{segment-based} (SB) and the \textit{element-based} (EB) schemes (Figure \ref{fig:integration_scheme}).  
The SB scheme builds a common partition intersecting the master and slave mesh. In this way, both basis functions are smooth over each element of the new grid and can be integrated exactly with an appropriate Gauss quadrature rule. However, finding and discretizing the geometric intersection of the master and slave grids can require a considerable implementation effort, especially in general 3D settings.
For this reason, the SB scheme is almost never used in general-purpose 3D codes.
By distinction, in the EB scheme integration is directly performed on the slave elements, avoiding the difficulty of  intersecting general meshes. The method is much easier to implement, but a quadrature error can arise since the master basis functions are not smooth in the integration domain. 

A common feature of the two approaches is the definition of the operator $\Pi$ as a geometric projection.
In order to compute $\Pi N_l^{\Gamma_1}$ on the quadrature points identified in the slave discretization,
each Gauss point must be projected onto the master side to evaluate the master basis function. This generally involves the solution of a non-linear problem, which can be performed by means of few Newton iterations. In the EB scheme, points projecting outside the master element are sorted out from integration.
For a more in-depth comparison between SB and EB scheme, the reader is referred to \cite{farah2015segment}, where the implementation and the accuracy of the two approaches have been investigated also in contact mechanics applications. 
The results there reported confirm that optimal convergence is preserved only with exact integration \cite{maday2002influence}. Nonetheless, in typical mesh tying problems, convergence of the EB scheme deviates from SB only for extremely fine meshes and well below the accuracy typically required by engineering applications. Hence, in practice the EB scheme usually proves to be a very effective and robust alternative to the SB approach.

Unfortunately, there are situations where the EB scheme can exhibit a severe loss of accuracy. This may occur at the boundary of the interfaces, whenever a slave elements is only partially included in the master support. Another critical scenario can be encountered in frictional contact applications. Both cases are discussed in detail in \cite{farah2015segment}. 
Obviously, increasing the number of Gauss quadrature points is the most straightforward way to mitigate the loss of accuracy, but the method may rapidly become inefficient because of the high computational cost required by the projection operations.

\section{A novel quadrature algorithm for mortar integrals} \label{sect3:RBFalgorithm}
The observations made in the previous section motivate the development of an alternative EB approach, for which the evaluation of the projection of the master basis functions at the integration points becomes a much cheaper operation.
In particular, we aim to avoid the projection of the Gauss quadrature points, which involves the solution of a set of non-linear problems.
The main idea consists in redefining $\Pi$ as an interpolation operator, instead of a projector. 
To this aim, we employ Radial Basis Functions (RBF) to compute a mesh-free interpolation of each master basis function, which can be easily evaluated at any point on the slave side. 
%
%
%
The RBF interpolation introduces an additional approximation in the EB approach. 
We will show that the possible loss of accuracy is generally small and can be kept under control, while maintaining the benefits of improved algorithmic efficiency as the number of integration points increases.

The proposed approach is denoted in the following as \textit{Radial-Basis} (RB) scheme. 

In this section, we briefly recall the fundamental properties of the RBF interpolation and show how it can be employed within the mortar integral. Then, we discuss the algorithmic implementation.

\subsection{RBF interpolation}

Consider a function $f: \Omega \subset \mathbb{R}^d \rightarrow \mathbb{R}$ and the set $\Xi = \{\boldsymbol \xi_m \}_{m=1}^{M}$ of interpolation points in $\mathbb{R}^d$, with $f_m = f(\boldsymbol{\xi}_m)$. The RBF interpolation defines a global interpolant $\Pi_f(\boldsymbol x)$ in $\Omega \subset \mathbb{R}^d$ as
\begin{align}
    \Pi_f(\boldsymbol x) = \sum_{m=1}^{M} \gamma_m^f \phi(\| \boldsymbol x - \boldsymbol{\xi}_m \| ,\varepsilon)\, ,
    \label{rbf_definition}
\end{align}
where $\phi(\| \boldsymbol x - \boldsymbol{\xi}_m \| ,\varepsilon)$ is the RBF, $\varepsilon$ is the so-called \textit{shape parameter}, and $\{ \gamma_m^f\}^{M}_{m=1}$ are the interpolation weights, computed imposing the interpolation condition
\begin{align}
    \Pi_f(\boldsymbol \xi_m) = f_m\, .
    \label{RBF_interpolationCondition}
\end{align}
By collecting $\{f_m\}_{m=1}^{M}$ and $\{ \gamma_m^f\}^{M}_{m=1}$ in the vectors $\boldsymbol f$ and $\boldsymbol \gamma_f$, equation \eqref{RBF_interpolationCondition} yields the linear system:
\begin{align}
    \Phi_{MM} \boldsymbol \gamma_f = \boldsymbol f
    \label{linsyst_RBF}
\end{align}
where $\Phi_{MM}[i,j] = \phi(\| \boldsymbol \xi_i - \boldsymbol \xi_j \| ,\varepsilon)$.
Possible available options for the RBF definition are reported in Table \ref{tabRBF}. Inverted Multiquadratic (IMQ) and Gaussian splines (GA) have a global support, thus producing a dense matrix $\Phi_{MM}$ for every choice of $\varepsilon$. By distinction, Wendland $C^2$ RBF has a local support, with $\Phi_{MM}$ having a zero entry for each pair of points whose distance is larger than $\varepsilon$.

\begin{table}
\centering
\caption{Examples of Radial Basis Functions $\phi(\| \boldsymbol x \| ,\varepsilon)$.}\label{tabRBF}%
\begin{tabular}{l c}
\toprule
Name &  $\phi(\| \boldsymbol x \| ,\varepsilon)$\\
\midrule
Gaussian splines (GA)    & $\exp({-\frac{\| \boldsymbol x \|^2}{
\varepsilon^2}})$   \\ [2pt]
Inverted multiquadratic (IMQ) & $(\| \boldsymbol x \|^2 + \varepsilon^2)^{-1/2}$   \\[2pt]
Wendland $C^2$ & $( 1 - \frac{\| \boldsymbol x \|}{\varepsilon})_{+}^4 ( 1 + 4 \frac{\| \boldsymbol x \|}{\varepsilon})$\\[1pt]
\bottomrule
\end{tabular}
\footnotetext[1]{$(x)_{+} = x \ \text{if} x > 0 \text{else} (x)_{+} = 0 \ \text{if} x \leq 0$}
\end{table}

Once the weights have been computed by solving \eqref{linsyst_RBF}, the interpolant can be evaluated on another set of $N$ points $\{\boldsymbol \zeta_n\}_{n=1}^{N}$:
\begin{align}
    \boldsymbol f_{\zeta} = \Phi_{NM} \boldsymbol \gamma_f \,.
    \label{mat_vec_interpolant}
\end{align}
Vector $\boldsymbol f_{\zeta}$ collects $\{\Pi_f(\boldsymbol \zeta_n)\}_{n=1}^{N}$, while $\Phi_{NM}[i,j] = \phi(\| \boldsymbol \zeta_i - \boldsymbol \xi_j \| ,\varepsilon)$.
In the following, the interpolant will be rescaled employing the procedure first introduced in \cite{deparis2014rescaled} for Wendland RBF, i.e.:
\begin{align}
    \hat\Pi_f(\boldsymbol x) = \frac{\Pi_f(\boldsymbol x)}{\Pi_1(\boldsymbol x)} = \frac{\sum_{m=1}^{M} \gamma_m^f \phi(\| \boldsymbol x - \boldsymbol{\xi}_m \| ,\varepsilon)}{\sum_{m=1}^{M} \gamma_{1,m} \phi(\| \boldsymbol x - \boldsymbol{\xi}_m \| ,\varepsilon)} \, ,
    \label{rescaling}
\end{align}
where $\Pi_1(\boldsymbol x)$ is the RBF interpolant of the constant function $f=1$ and $\{ \gamma_{1,m}\}^{M}_{m=1}$ the corresponding weights. Rescaling \eqref{rescaling} requires solving one additional linear system with the matrix $\Phi_{MM}$, but improves accuracy and yields exact reproduction of constant functions. Linear convergence of the rescaled RBF has been demonstrated in \cite{de2020convergence} up to a certain conjecture. In this work, the rescaling \eqref{rescaling} will be applied also to other families of RBF.

\begin{remark}
RBF interpolants may exactly reproduce polynomials of any degree $k$ by enriching \eqref{rbf_definition} with a polynomial term $p(\boldsymbol x) \in \mathbb{P}_k(\mathbb{R}^d)$. To ensure that $\Phi_{MM}$ remains non singular, the set $\Xi$ must be $k$-unisolvent, i.e., there exists a unique polynomial in $\mathbb{P}_k(\mathbb{R}^d)$ of lowest possible degree that interpolates the data in $\Xi$. 
However, in our case this is not a viable option. For instance, when $\Omega \subset \mathbb R^2$, the interpolation points will be located on straight segments, thus representing a non-unisolvent set for polynomials in $\mathbb{P}_1(\mathbb{R}^2)$. Equivalently, when $\Omega \subset \mathbb R^3$, interpolation points in the master elements are coplanar and not-unisolvent for $\mathbb{P}_1(\mathbb{R}^3)$.
\end{remark}

\subsection{Evaluating mortar integrals with the RB approach} \label{sect:eval_int_RBF}

Let us consider an element lying on the master side $\Gamma_1$ with a set $\{\boldsymbol \xi_m\}_{m=1}^{M}$ of interpolation points, and an element lying on the slave side $\Gamma_2$ with $N_G$ Gauss points $\{\boldsymbol \zeta_g\}_{g=1}^{N_G}$.
First, we interpolate the nodal basis functions $N^{\Gamma_1}_j$ on a master element computing the RBF weights.
If we denote by $N_{\Gamma_1}$ the local matrix of basis functions such that $N_{\Gamma_1}[i,j] = N_{j}^{\Gamma_1}(\boldsymbol{\xi}_i)$, the matrix of RBF weights reads:
\begin{equation}
    W = \Phi_{MM}^{-1} N_{\Gamma_1} . 
\label{master_interpolation}
\end{equation}
The RBF weights of the interpolant of the unitary function are
\begin{equation}
    \boldsymbol w_1 = \Phi_{MM}^{-1} \boldsymbol 1_{M} \, ,
\label{unitary_weights}
\end{equation}
where $\boldsymbol 1_M$ is the vector of $M$ ones. The scaled interpolant is then evaluated at the Gauss quadrature points, obtaining:
\begin{align}
    N_{\Gamma_2}^{\Pi} &= R^{-1}\Phi_{NM} W
    \label{Interpolant_eval}
\end{align}
with $R = diag(\Phi_{NM} \boldsymbol w_1$). The resulting matrix $N_{\Gamma_2}^{\Pi}$ collects the values of the scaled RBF interpolant of the master basis function at the Gauss points, i.e., $N_{\Gamma_2}^{\Pi}[i,j] = \Pi N_{j}^{\Gamma_1}(\boldsymbol{\zeta}_i)$. The rescaling operation through matrix $R$ guarantees that $N_{\Gamma_2}^{\Pi}$ preserves the partition of unity property, 
hence the method can be consistent (see Remark \ref{rem:unity}). 

Once the projections of the master basis functions at the Gauss points are available, an EB integration is performed on the slave elements to compute the mortar matrices.
Similarly as before, let us denote by $N_{\Gamma_2}$ and $\Psi$ the matrices collecting the values of the slave and the Lagrange multiplier basis functions on the Gauss points, $N_{\Gamma_2}[i,j]=N_j^{\Gamma_2}(\boldsymbol{\zeta}_i)$ and $\Psi[i,j]=\psi_j(\boldsymbol{\zeta}_i)$, and by $W_g$ the $N_G\times N_G$ diagonal matrix with entries equal to the product of Gauss weight by the Jacobian of the isoparametric map computed at each quadrature point $\boldsymbol{\zeta}_i$. 
From such matrices, we drop the rows corresponding to the Gauss points selected in the slave element that fall outside the master element, thus reducing the row dimension of all matrices to $N_G^* \le N_G$. 
Details on how to determine the $N_G^*$ active quadrature points are provided in \ref{implementation_details}.
Finally, the mortar matrices read:
\begin{align}
    D = 
    \Psi^T W_g N_{\Gamma_2} \, , \label{D_eb} \\
    S = 
    \Psi^T W_g N_{\Gamma_2}^{\Pi} \, . \label{M_eb}
\end{align}
%
For the practical implementation of the above procedure, we need to define:
\begin{enumerate}
    \item the type of RBF;
    \item the value of the shape parameter $\varepsilon$;
    \item the number $M$ and location of interpolation points in the master element;
    \item the number $N_G$ of Gauss points in the slave element.
\end{enumerate}

\begin{figure}
    \centering
    \includegraphics[width=13.5cm]{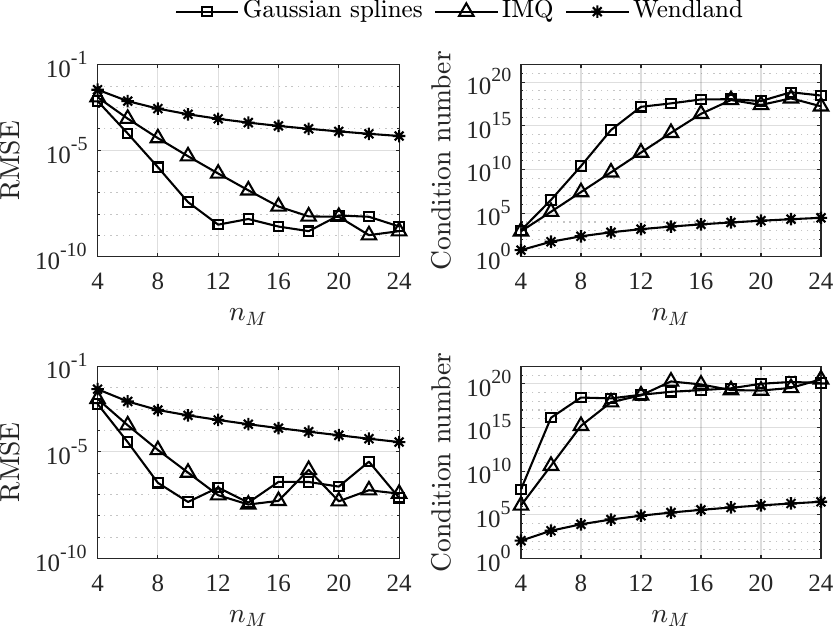}
    \caption{Accuracy and condition number for interpolating the shape functions of 1D quadratic elements (top) and 2D serendipity elements (bottom).}
    \label{figBasisF_interp}
\end{figure}

\begin{figure}
    \centering
    \includegraphics[width=0.9\linewidth]{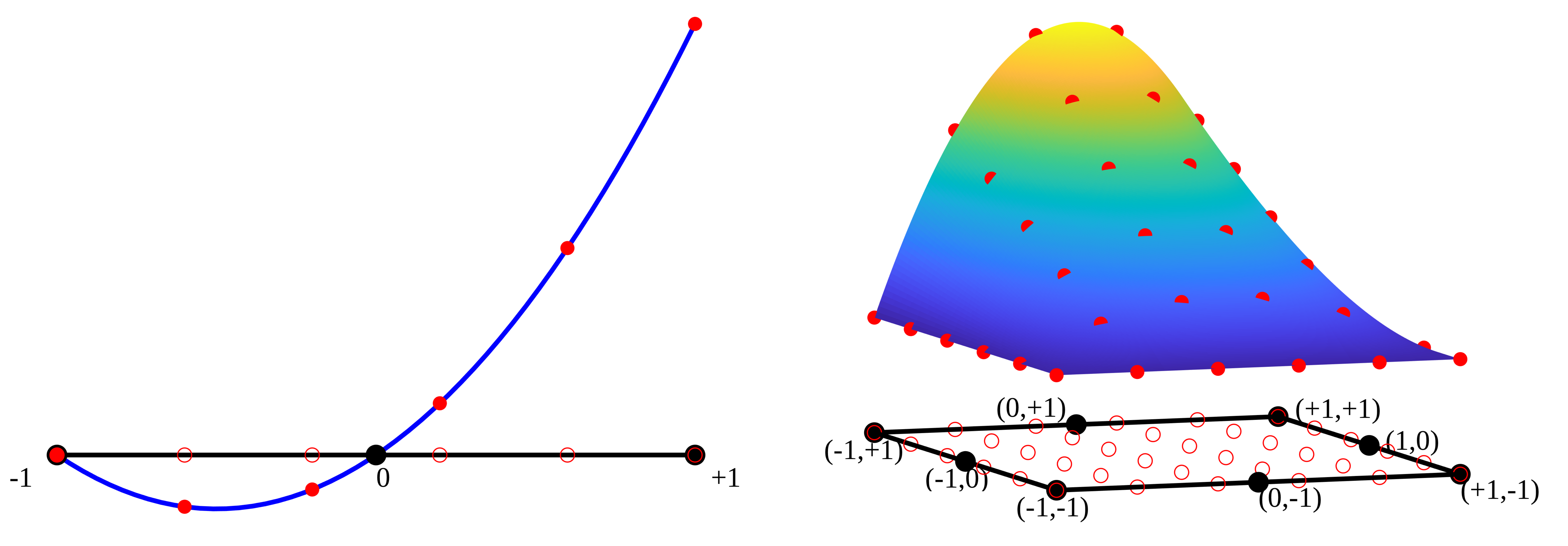}
    \caption{Quadratic reference elements in 1D (3 nodes) and 2D (8 nodes) for RBF comparison. Red dots represent interpolation points in a regular grid. In this figure $n_M = 6$.}
    \label{Ref_elements}
\end{figure}

The type of RBF can significantly affect the resulting accuracy. Of course, we aim at obtaining the required accuracy keeping $M$ as small as possible. The RBF classes with a global support are a natural choice, since $\Phi_{MM}$ is small and memory usage, as well as matrix inversion, is not a concern. Figure~\ref{figBasisF_interp} compares the \textit{root-mean-square-error} (RMSE),
\begin{align*}
RMSE = \bigg[ \frac{1}{N} \sum_{n=1}^N (\hat\Pi_f(\boldsymbol \zeta_n) - f(\boldsymbol \zeta_n))^2 \bigg]^{1/2},
\end{align*}
and the condition number of $\Phi_{MM}$, when interpolating 1D quadratic and 2D biquadratic basis functions in the reference finite elements using the RBF options listed in Table \ref{tabRBF}. The set $\Xi$ of interpolation points is taken by a uniform subdivision with $n_M$ points on each edge (Figure \ref{Ref_elements}), implying $M = n_M$ for 1D elements, $M = n_M\cdot(n_M+1)/2$ for triangles, and $M=n_M^2$ for quadrilaterals, while the points $\boldsymbol \zeta_n$ used for the RMSE evaluation are taken from a set of $N=40$ Halton points. The shape parameter $\varepsilon$ is set equal to the circumdiameter of the element.
As shown in Figure \ref{figBasisF_interp}, GA and IMQ can provide a very accurate interpolation of the polynomial basis functions with a limited number of points. However, matrix $\Phi_{MM}$ gets rapidly ill-conditioned as $\varepsilon/M$ decreases. As a consequence, the interpolation becomes unstable and fails to converge further. To overcome this issue, known as \textit{trade-off principle}, research has focused on finding ``optimal'' values for $\varepsilon$ \cite{rippa1999algorithm}, as well as developing alternative stable RBF algorithms when $\varepsilon/M \rightarrow 0$ \cite{fornberg2011stable,fornberg2004stable}. Anyway, such an effort can be avoided in our method, since a satisfactory accuracy is usually already achieved in the stable regime. 
Unlike globally supported kernels, the class of Wendland RBF remains stable even for small $\varepsilon/M$, so they might be considered a robust alternative when a smaller accuracy suffices. 

\begin{figure}
    \centering
    \includegraphics[height=7cm]{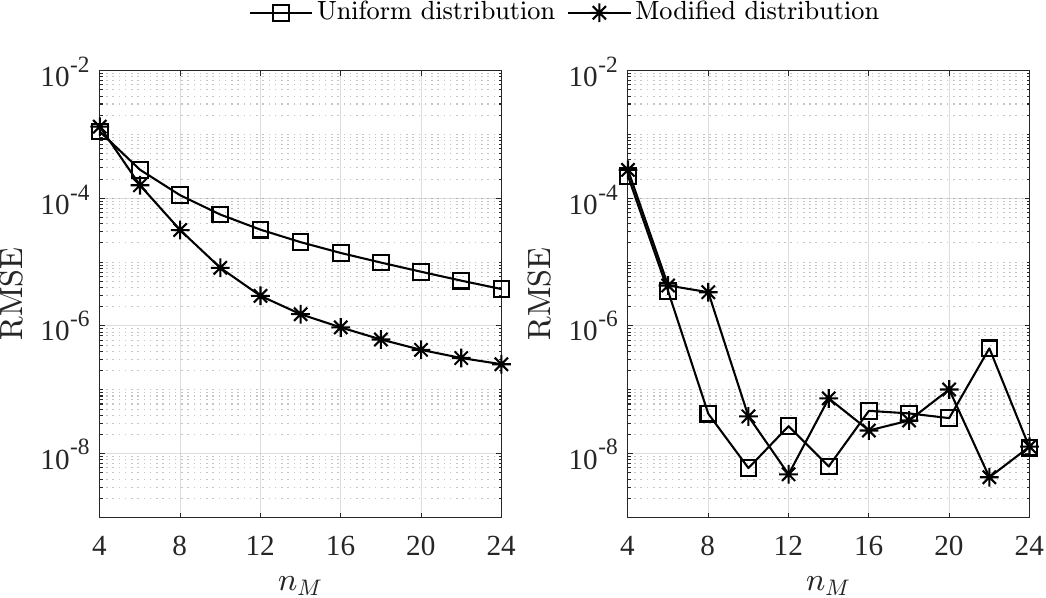}
    \caption{Influence of the dataset location on the interpolation error: Wendland RBFs at left, Gaussian splines at right.}
    \label{RBF_dataSet}
\end{figure}

Regarding the location of the interpolation points, the previous test is repeated with GA and Wendland RBF comparing the outcome obtained from equally spaced points ($\Xi_{uni}$, Figure \ref{Ref_elements}) with the modified distribution  
\begin{align*}
    \Xi_{mod} = \bigg\{ \boldsymbol\xi_m^{mod} = \sin{\bigg(\frac{\pi}{2}\boldsymbol\xi_m\bigg)} \ \vert \ \boldsymbol\xi_m \in \Xi_{uni} \bigg\}.
\end{align*}
The results are shown in Figure \ref{RBF_dataSet}. GA is barely affected by the dataset choice, while Wendland RBF improves by about one order of magnitude with $\Xi_{mod}$. This can be explained by the fact that points in $\Xi_{mod}$ are denser near the boundaries, where Wendland RBFs tend to exhibits larger errors.

GA not only provides a superior accuracy in the stable regime, but also ensures a greater interpolation robustness in the presence of gaps between the master and slave interfaces. This is motivated by the result that follows.

\begin{proposition}
\emph{
Let $\Xi=\{\boldsymbol \xi\}_{m=1}^{M}$ be a set of collinear (in 1D) or coplanar (in 2D) points defined on a master element with unit normal $\boldsymbol{n}$. The rescaled interpolant \eqref{rescaling} computed with Gaussian splines yields:}
\begin{align}
    \hat\Pi_f(\boldsymbol{x}) =  \hat\Pi_f(\boldsymbol{x} + s\boldsymbol{n}), \ \ \ \ \ \ \forall \, s \in \mathbb{R}, 
    \label{geom_flexibility}
\end{align}
\emph{for any point $\boldsymbol x$ lying on the master element.}
\end{proposition}

\begin{proof}
The GA exponential form yields
\begin{eqnarray}
    \phi(\boldsymbol{x} + s\boldsymbol{n} - \boldsymbol \xi_m) &=& \phi(\boldsymbol{x} - \boldsymbol \xi_m) \, \phi(s\boldsymbol{n}) \, \phi(\sqrt{2s(\boldsymbol \xi_m - \boldsymbol x)^T \boldsymbol n}) \nonumber \\
    &=& \phi_{m}^{(1)}(\boldsymbol x)  \phi(s\boldsymbol{n}) \phi_{m}^{(2)}(\boldsymbol x).
    \label{eq:prop_GA}
\end{eqnarray}
The term $\phi(s\boldsymbol{n})$ is a constant in both the numerator and the denominator of equation \eqref{rescaling}, hence it cancels out after rescaling. 
Moreover, we notice that $(\boldsymbol \xi_m - \boldsymbol x)^T \boldsymbol n = 0$ for every point $\boldsymbol x$ of the master element, since $\boldsymbol \xi_m - \boldsymbol x$ lies in a plane with normal $\boldsymbol n$. This observation implies $\phi_m^{(2)}(\boldsymbol x) = 1$, thus concluding the proof.
\end{proof}

Equation \eqref{geom_flexibility} implies that the basis function interpolant is independent of any shift of the slave element along the master normal. In this way, the rescaled interpolation builds a projection onto the master surface and guarantees accuracy also in geometrically non-conforming settings.
For these reasons, in the sequel we will consider GA as the default RBF selection.

Finally, we address the influence of the shape parameter $\varepsilon$. 
Two possible choices are compared: (i) fixed $\varepsilon$ proportional to the circumdiameter of the master element, indicated with $h_M$; (ii) $\varepsilon$ changing with $M$ and proportional to the fill distance $h_{\Xi}$ of the dataset, defined as
\begin{align*}
    h_{\Xi} = \sup\limits_{\boldsymbol x \in \Omega} \, \min\limits_{\boldsymbol \xi_m \in \, \Xi} \| \boldsymbol x - \boldsymbol \xi_m \|_2.
\end{align*}
\begin{figure}
    \centering
    \includegraphics[height=6.5cm]{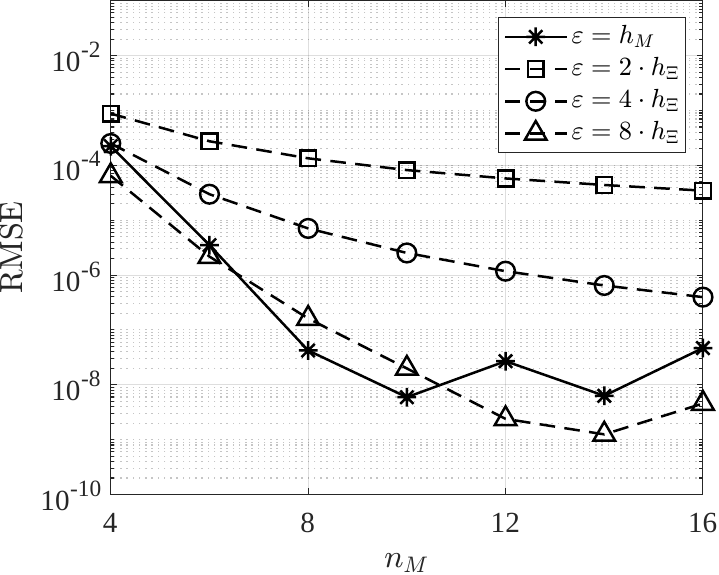}
    \caption{Influence of the shape parameter $\varepsilon$ on the interpolation error using GA.}
    \label{test_radius}
\end{figure}
Figure \ref{test_radius} shows the outcomes of these two choices for the interpolation of a quadratic basis function in 3D. Setting $\varepsilon=h_M$ yields a good trade-off between accuracy and stability when $n_M$ is small. When $\varepsilon$ varies with the fill distance, the accuracy is significantly affected by the ratio $\varepsilon/h_{\Xi}$. 
Better results are obtained increasing $\varepsilon$, but instabilities occur already for $\varepsilon/h_{\Xi} = 8$. 
Therefore, from now on we will stick to the simple choice $\varepsilon=h_M$. 

Based on these results, the proposed RB integration scheme will be configured as follows:
\begin{enumerate}
    \item Gaussian splines as RBF family;
    \item uniformly distributed interpolation points in the master element reference space;
    \item fixed shape parameter equal to the circumdiameter of the master element.
\end{enumerate}
Therefore, in our algorithm the only real parameters are $n_M$ (hence $M$) and $N_G$. 
Based on our numerical experiments, the restriction $n_M\leq10$ is suggested to avoid the occurrence of instabilities. 

\subsection{Implementation details} \label{implementation_details}
Algorithm \ref{alg:RB_scheme} summarizes the straightforward implementation of the proposed RB scheme, which can be extended from 2D to 3D without any difficulty. As in any mortar implementation, a contact search algorithm is needed to find each pair of master/slave elements sharing support (line 8 of Algorithm \ref{alg:RB_scheme}). State-of-the-art strategies are presented in \cite{klosowski1998efficient,yang2008contact}. The rest of the procedure involves standard algebraic computations, already described in section \ref{sect:eval_int_RBF}. In order to detect if a Gauss point on a slave element is located within the support of a corresponding master element, we use a \textit{support detection} procedure (line 13 of Algorithm \ref{alg:RB_scheme}). In the EB scheme, a Gauss point belongs to the master support if its projection lies inside the reference element domain. 
In the RB scheme, support detection is carried out in a similar way. 
If a nodal basis functions $N^{\Gamma_1}(\boldsymbol{x})$ is non-negative within a master element, then also $\Pi N^{\Gamma_1}(\boldsymbol{x}_g)$ must be non-negative on those slave Gauss points projecting onto the master support. Points not fulfilling this condition are removed from EB integration. 
This procedure is robust and may be plugged also in standard EB algorithms, in order to avoid performing unnecessary projections outside the master support.


\begin{remark}
In the case of basis functions that are not always positive in the element, auxiliary positive functions need to be interpolated solely for the purpose of support detection. For quadrilateral elements of any order, such auxiliary interpolation is required only for basis functions of opposite corner nodes. 
\end{remark}

\begin{algorithm}
\caption{RB quadrature scheme.}
    \begin{algorithmic}[1]
    \Procedure{RB\_integration}{$n_M$,$N_G$}
    \For{$l = 1,\dots ,\vert \mathcal{N}_{\Gamma_1} \vert$ }    \Comment{Master interpolation loop}
    \State Locate interpolation points in the master element
    \State Compute $\Phi_{MM}$ and $N_{\Gamma_1}$
    \State Compute $W$ and $\boldsymbol{w}_{1}$ \Comment{equations \eqref{master_interpolation}-\eqref{unitary_weights}}
    \EndFor
    \For{$i = 1,\dots ,\vert \mathcal{N}_{\Gamma_2} \vert$ } \Comment{Integration on the slave interface}
    \State Find connected master elements with contact search algorithm
    \State Define $N_G$ Gauss points and compute multiplier basis functions
    \ForAll{connected master elements}
    \State Compute $\Phi_{NM}$
    \State Compute interpolant $N_{\Gamma_2}^{\Pi}$ \Comment{equation \eqref{Interpolant_eval}}
    \State Sort out Gauss points with support detection
    \State Compute local contribution to $D$ \Comment{equation \eqref{D_eb}}
    \State Compute local contribution to $S$ \Comment{equation \eqref{M_eb}}
    \EndFor
    \EndFor
    \EndProcedure
    \end{algorithmic}
    \label{alg:RB_scheme}
\end{algorithm}

\begin{remark}
Notice that in Algorithm \ref{alg:RB_scheme} no interpolation is performed to compute the multiplier basis functions, since they are already available on the slave element. Due to the interpolation error introduced in $S$, if we directly apply Algorithm \ref{alg:RB_scheme} to a conforming interface, $D$ and $S$ could not exactly coincide, hence $E\neq I$, with $I$ the identity matrix. 
Anyway, this observation has no practical implications, since in the conforming case the mortar algorithm is not needed.
\end{remark}

\section{Numerical results} \label{sect4:NumericalResults}
In this section the proposed RB scheme is validated and tested in a set of numerical experiments. First, we analyze accuracy and convergence of the RB scheme as compared to the standard EB scheme in 2D and 3D mortar interface problems where a known function defined on the master side is transferred on the slave side. Then, 2D and 3D Poisson problems are solved subdividing the domain into non-conforming discretizations and applying the proposed mortar algorithm at the common interface.
Finally, a realistic application of the proposed method for the solution to a fluid-structure interaction problem is presented.

\subsection{Accuracy and convergence}
\begin{figure}
    \centering
    \includegraphics[width=13cm]{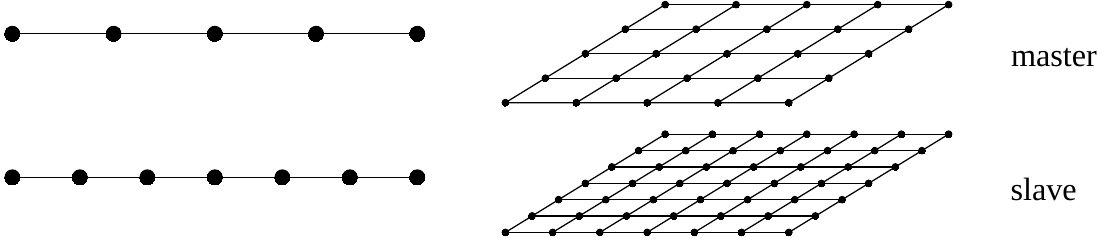}
    \caption{Examples of discretization of 1D and 2D geometrically conforming interfaces.}
    \label{Res_interp_mesh}
\end{figure}

\begin{figure}
    \centering
    \includegraphics[height=7.5cm]{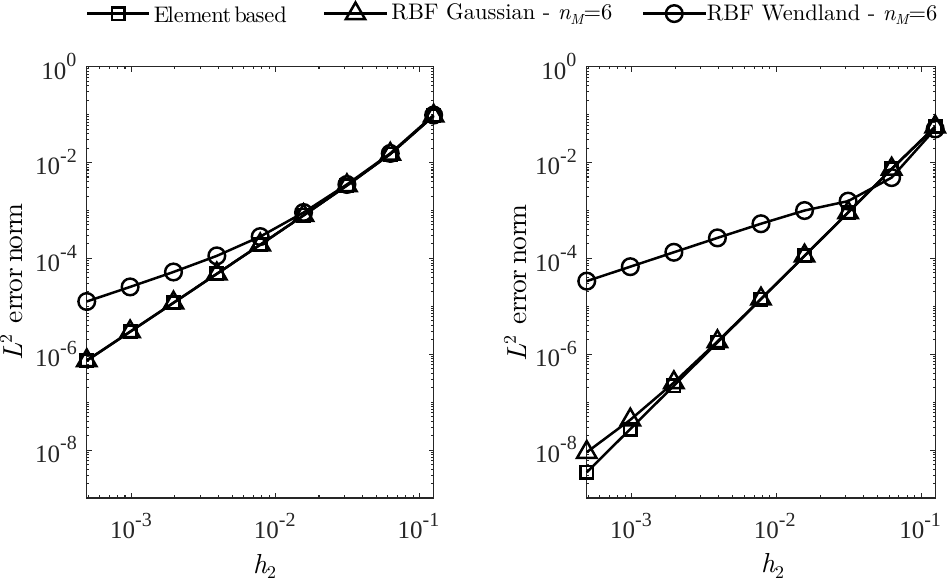}
    \caption{Convergence profiles for the 1D interpolation test: linear elements (left, 2 Gauss points) and quadratic elements (right, 3 Gauss points).}
    \label{interp_1D}
\end{figure}

First, we analyze the outcome obtained on a geometrically conforming interface (linear or planar) discretized with different grids on the master and slave side (Figure \ref{Res_interp_mesh}). Consider the 1D interval $[-1,+1] \subset \mathbb R$ discretized by either linear or quadratic elements. Let $h_1$ and $h_2$ denote the characteristic mesh sizes for master and slave interfaces, respectively, with $h_2/h_1=2/3$ kept constant during the convergence test (Figure \ref{Res_interp_mesh}, leftmost frame).
The analytical function $f(x) = \sin(4x)+x^2$ is defined on the master side and transferred to the slave side by using the mortar interpolation \eqref{discrete_mortar_constr_alg}. 
The integration is performed with the minimum number of Gauss points, i.e., 2 in the linear case, 3 in the quadratic case.
The $L^2$ error norm is then computed on the slave side while progressively reducing $h_2$ and $h_1$. The resulting convergence profiles are shown in Figure \ref{interp_1D} using both GA and Wendland RBF with $M=n_M=6$ equally distributed interpolation points for each master element. 
As expected, the accuracy of the RBF interpolation controls the convergence rate of the interface problem.
The RB scheme with GA produces $L^2$ errors that are virtually identical to those of the standard EB scheme for both linear and quadratic element discretizations. This accuracy is guaranteed at a very small cost, since only 6 interpolation points are necessary. In contrast, Wendland RBFs diverge from the optimal error profile for the same number of interpolation points because of its lower accuracy (Figure \ref{figBasisF_interp}), following a linear convergence \cite{de2020convergence}. 


\begin{figure}
    \centering
    \includegraphics[height=7.5cm]{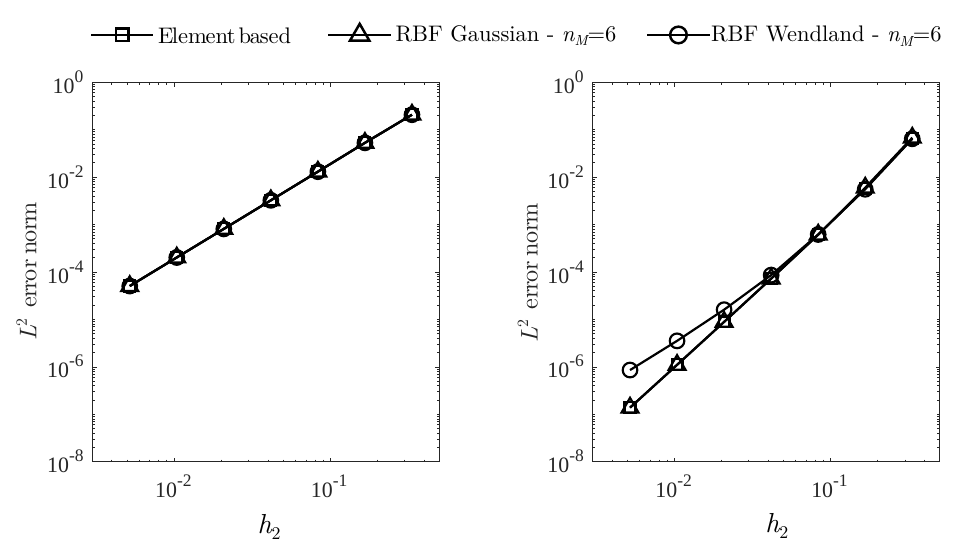}
    \caption{Convergence profiles for the 2D interpolation test: bilinear elements (4 Gauss points) and biquadratic serendipity elements (9 Gauss points).}
    \label{fig_interp2D}
\end{figure}

The same observations hold also in the 2D square interface $[-1,+1] \times [-1,+1] \subset \mathbb R^2$ of Figure \ref{Res_interp_mesh} (rightmost frame). The same assumption $h_2/h_1=2/3$ as for the 1D case is used. 
In this case, the integration is carried out with 4 and 9 Gauss points in the bilinear and biquadratic element, respectively.
The convergence profiles of Figure \ref{fig_interp2D} are obtained considering the interpolation of $f(x,y) = \sin(4x)\cos(4y)$ from the master surface by using bilinear and biquadratic quadrilateral elements. Here, only GA has been considered as RBF class, while the number of interpolation points has been modified from $n_M=4$ to $n_M=6$. Recall that for quadrilateral elements $M = n_M^2$. For a bilinear element discretization, both choices reproduce the outcome obtained with the EB scheme, while  only for a biquadratic element discretization with $n_M=4$ we can notice a slight deviation from the EB convergence for fine meshes.

\begin{figure}
\centering
\subfloat[]{\includegraphics[width = 7.2cm]{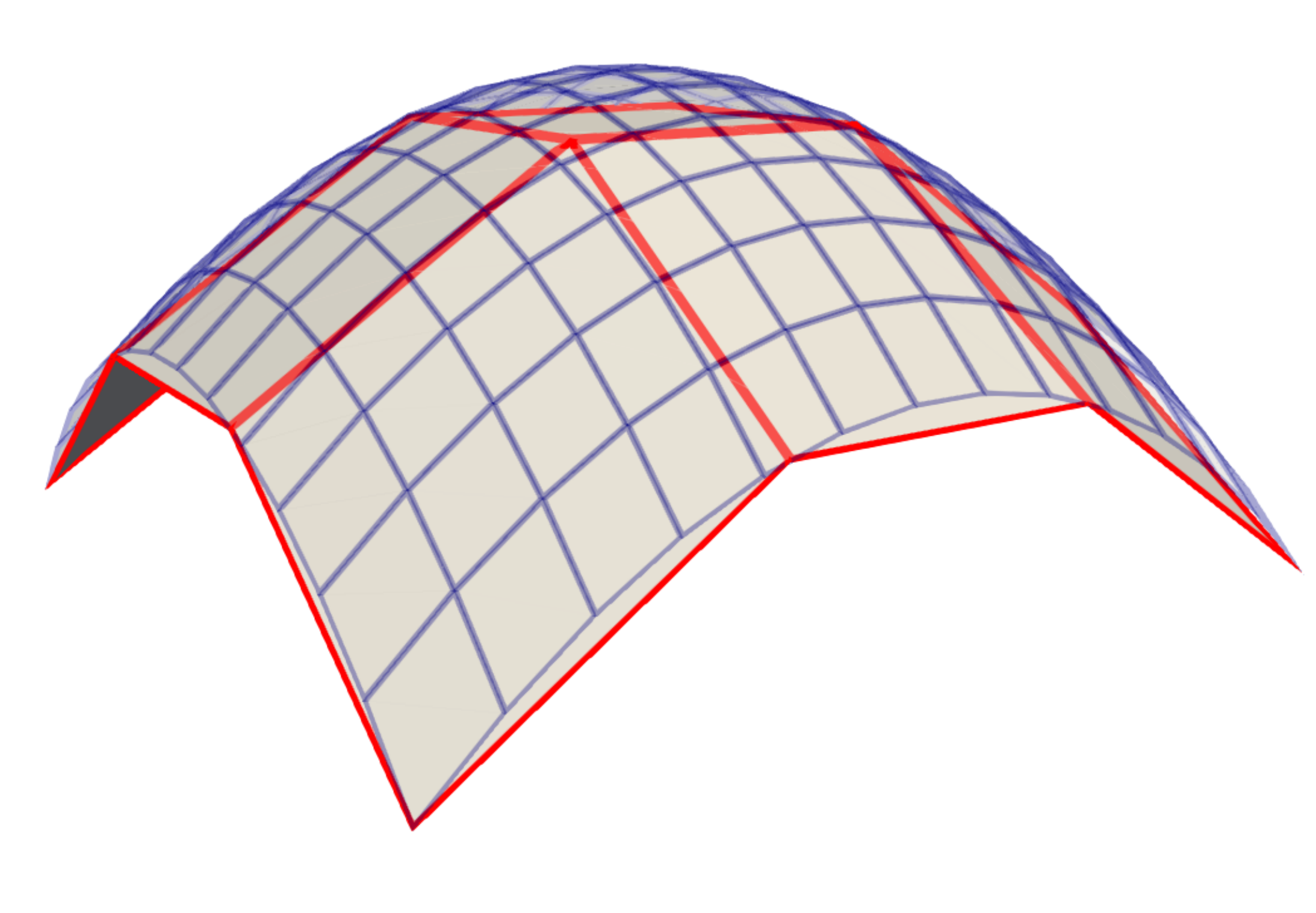}\label{fig:non_conf_meshes}}
\hspace{0.2cm}
\subfloat[]{\includegraphics[width = 6cm]{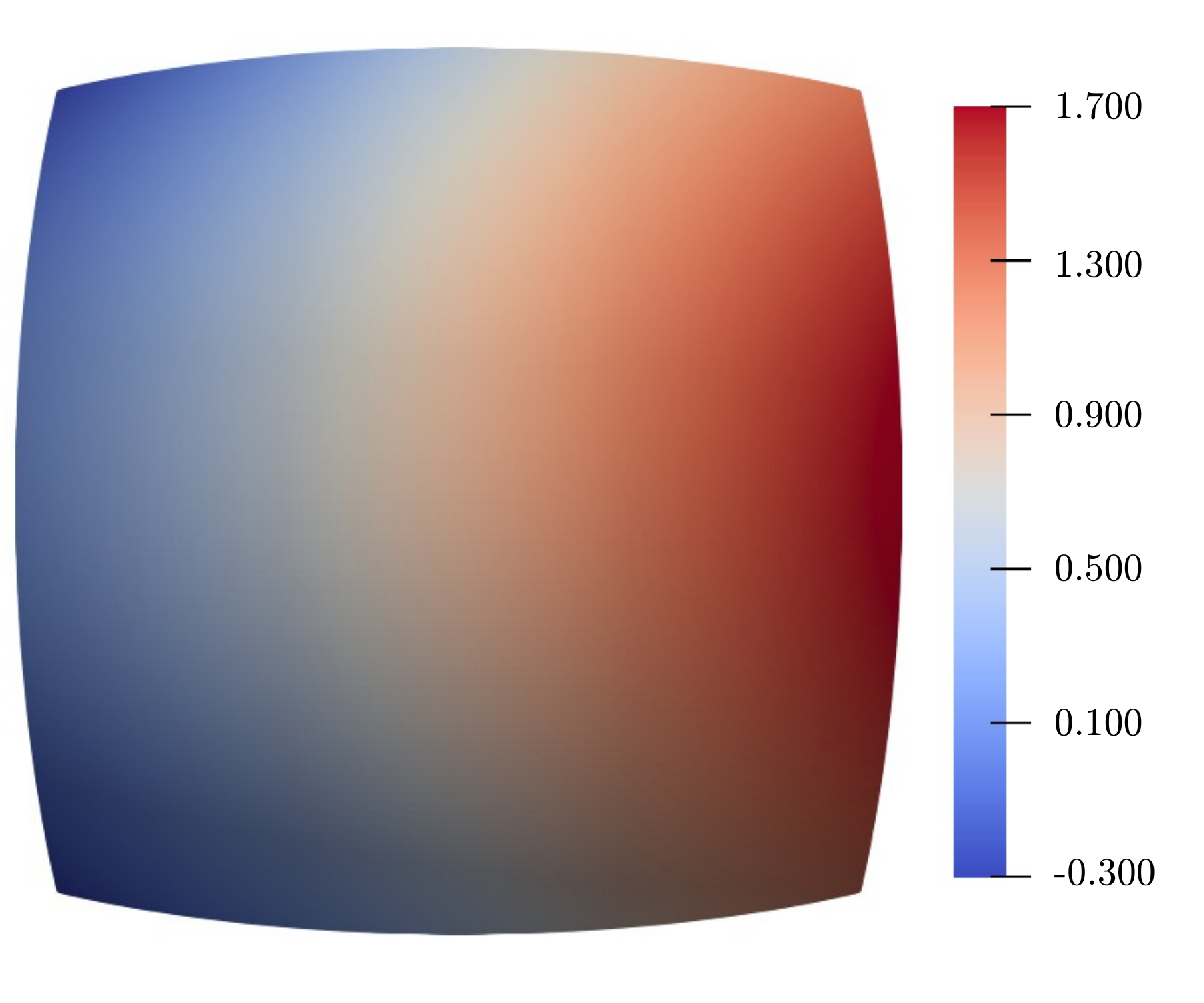}\label{fig:non_conf_anal}} 
\caption{Mortar interpolation on non-conforming interfaces: (a) discretization with a fine (blue edges) and a coarse grid (red edges); (b) contour plot of the function $f(x,y) = \sin(x)+\cos(y)$ defined on the master side.}
\end{figure}

\begin{figure}
    \centering
    \includegraphics[width=0.9\linewidth]{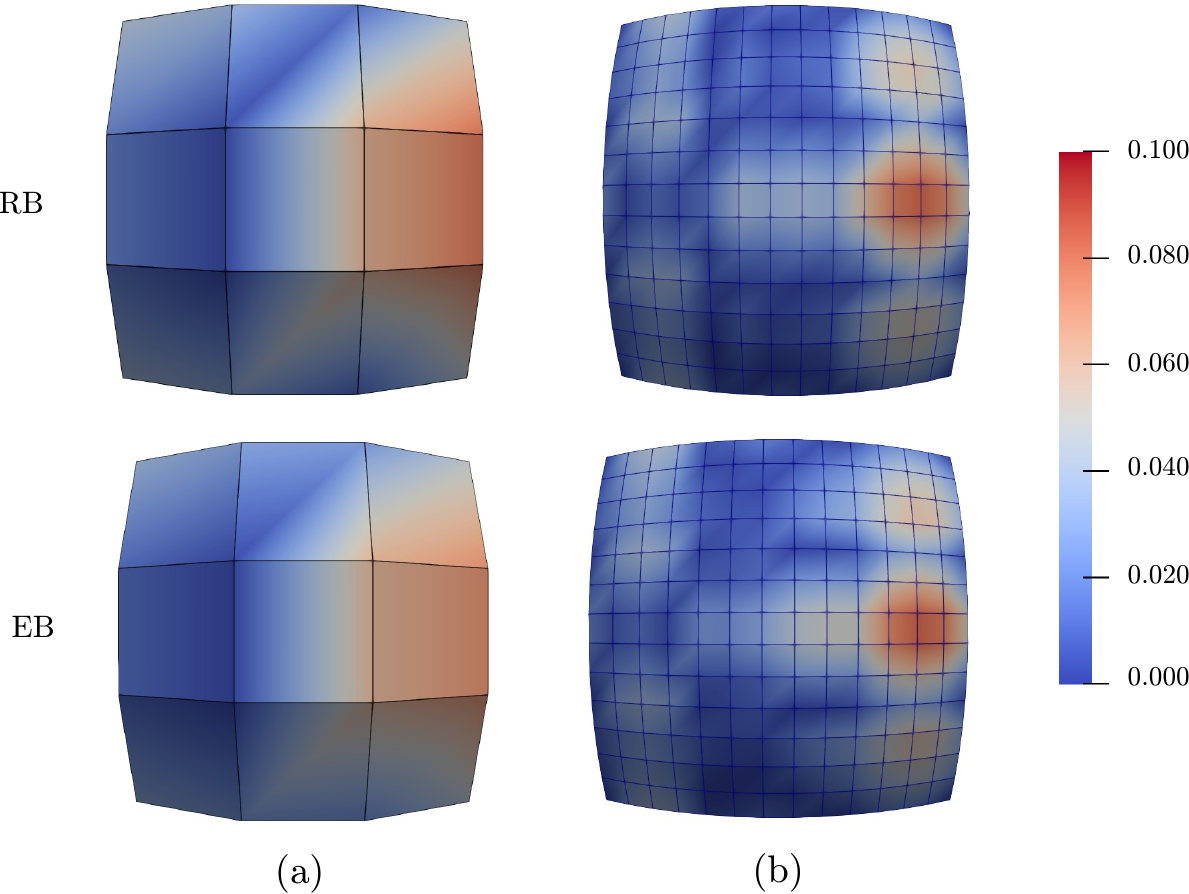}
    \caption{Mortar interpolation on non-conforming interfaces: absolute error on the slave side. (a) Coarse grid as slave domain, RB scheme on top ($L^2$ error = 0.0887), EB scheme on bottom ($L^2$ error = 0.0840); (b) fine grid on the slave domain, RB scheme on top ($L^2$ error = 0.0935), EB scheme on bottom ($L^2$ error = 0.0941).}
    \label{fig:NonConf_results}
\end{figure}

The 2D interface is now distorted to obtain the geometrically non-conforming interfaces shown in Figure \ref{fig:non_conf_meshes}. 
The function $f(x,y)=\sin(x)+\cos(y)$ (Figure \ref{fig:non_conf_anal}) is defined on the master side and then transferred on the slave side by using equation \eqref{discrete_mortar_constr_alg}.
Here, we compare the outcomes of the standard EB scheme with those obtained from the RB scheme with GA and $n_M = 4$ for the computation of the mortar operator $E$.
Figure \ref{fig:NonConf_results} shows the absolute error of the function $f(x,y)$ as reproduced on the slave side, set as either the coarse (Figure \ref{fig:NonConf_results}a) or the fine grid (Figure \ref{fig:NonConf_results}b).
Unsurprisingly, larger errors occur where master and slave surfaces exhibits larger gaps.
However, despite the challenging geometrical setting, the $L^2$ error norm is very close to the one obtained with EB scheme, so that the point error contours of Figure \ref{fig:NonConf_results} are almost indistinguishable. 
Notice that a better accuracy is achieved when the master side has the coarsest discretization. This result can be explained by the fact that gaps become less effective on the interpolation results, since they are smaller compared to the area of each master element.

\subsection{Poisson problems}
\begin{figure}
    \centering
    \includegraphics[width=\linewidth]{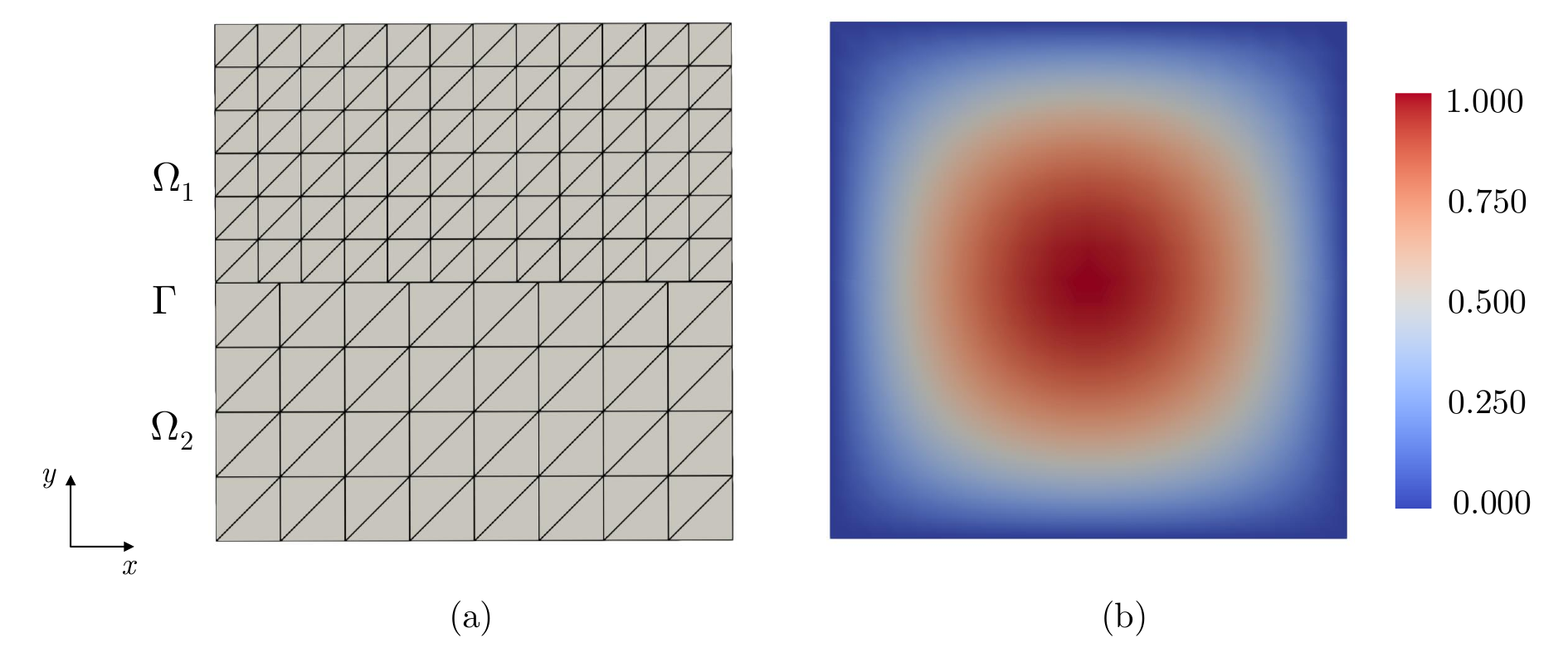}
    \caption{2D Poisson problem: problem setting. Example of discretized domain (a) and contour plot of the manufactured solution (b).}
    \label{fig:Poisson2D_setting}
\end{figure}

\begin{figure}
\centering
\subfloat[]{\includegraphics[width = 0.45\linewidth]{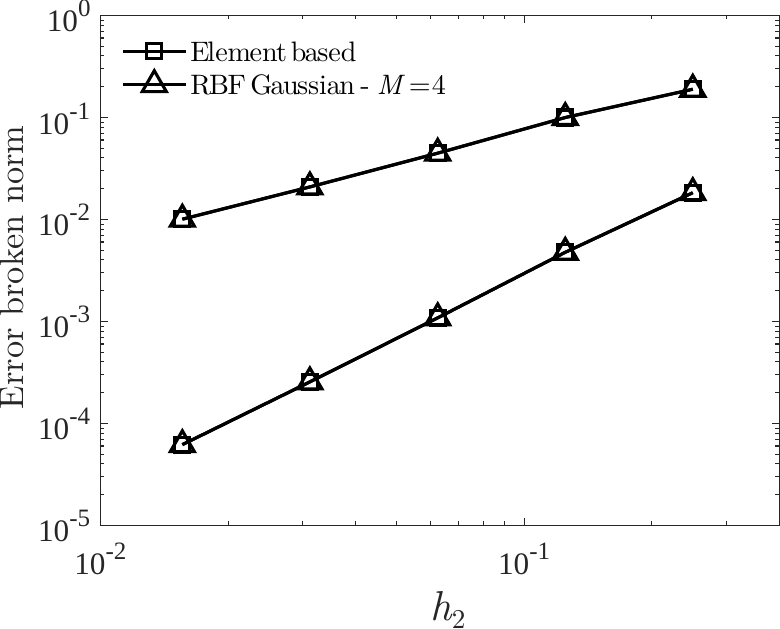}\label{Poisson2D_conv}}
\hspace{0.3cm}
\subfloat[]{\includegraphics[width = 0.45\linewidth]{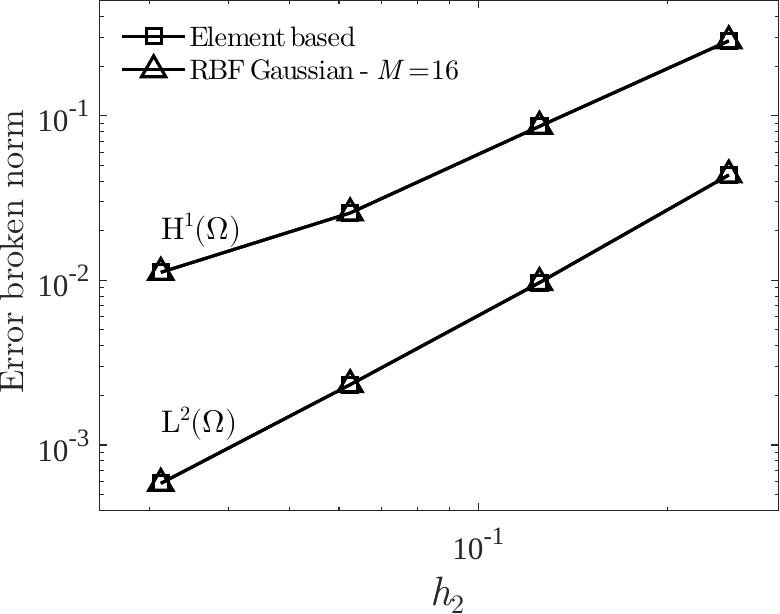}\label{Poisson3D_conv}} 
\caption{Convergence of Poisson problems in $L^2$ and $H^1$ broken norms: (a) 2D and (b) 3D case.}
\label{fig:Poisson_conv}
\end{figure}

\begin{figure}
    \centering
    \includegraphics[width=0.9\linewidth]{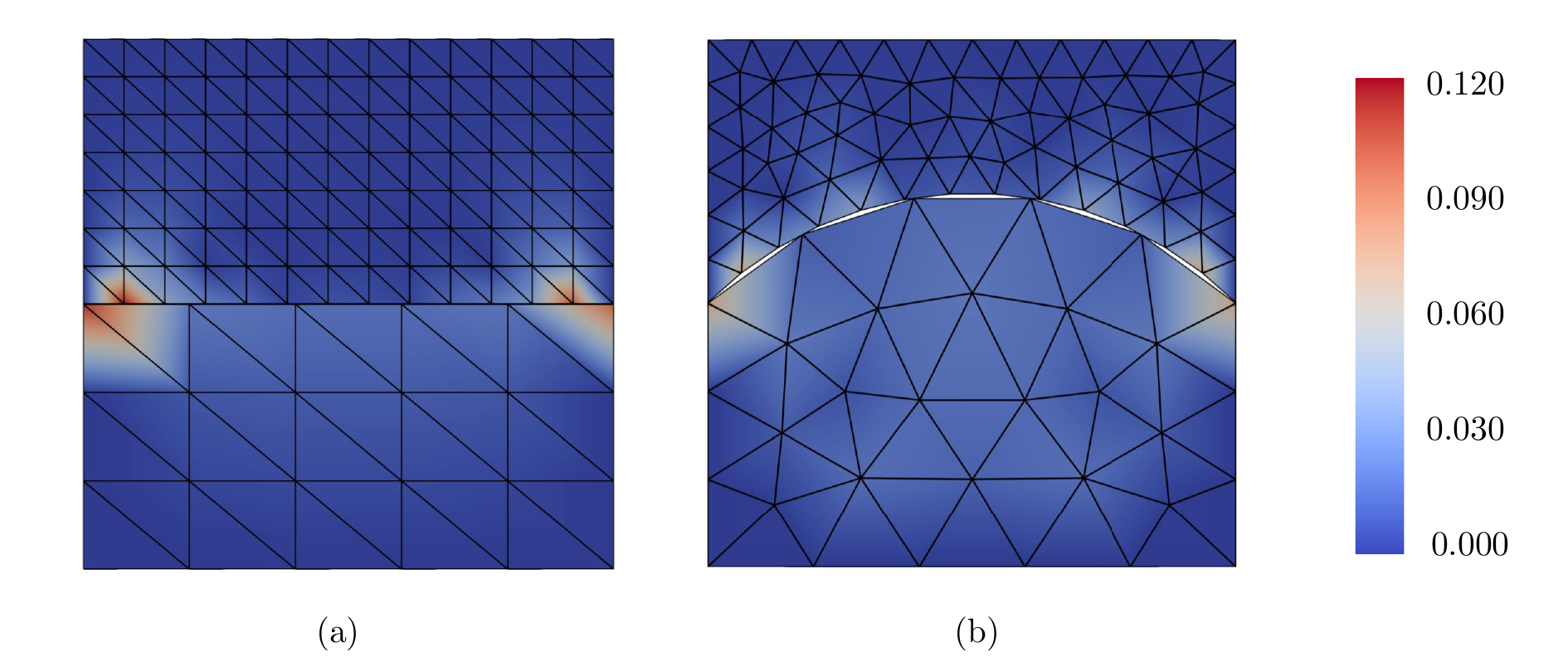}
    \caption{2D Poisson problem: absolute error of the solution with (a) a flat interface, and (b) a curved interface.}
    \label{fig:Pois2D_pointError}
\end{figure}

We consider the solution of the Poisson problem
\begin{equation*}
    -\Delta u = f \ \ \ \ \ \ \ \ \ \ \text{in} \ \Omega \subset \mathbb{R}^d, \, d=2,3
\end{equation*}
with forcing term $f$ and boundary conditions chosen to satisfy a certain manufactured solution.
A 2D setting is first considered. The domain $\Omega = [0,1] \times [0,1]$ is decomposed into $\Omega_1 = [0,1] \times [0.5,1]$ and $\Omega_2 = [0,1] \times [0,0.5]$. Corresponding partitions $\mathcal{T}_1$ and $\mathcal{T}_2$ are made of simplices, with $h_1/h_2 = 2/3$ at each uniform refinement (Figure \ref{fig:Poisson2D_setting}a). The manufactured solution $u(x,y) = 16xy(1-x)(1-y)$ is imposed (Figure \ref{fig:Poisson2D_setting}b) to compute $f$ and the corresponding essential boundary conditions. The error for evaluating the convergence of the discrete solution to the manufactured one is measured by the broken norm:
\begin{align*}
   \vertiii{v}_{H^r(\Omega)} = \sqrt{\sum_{k=1}^2 \| v_k \|_{H^r(\Omega_k)}}
\end{align*}
The convergence profiles reported in Figure \ref{Poisson2D_conv} confirm the consistency of the proposed mortar algorithm and show again that the accuracy obtained with the RB scheme is in practice indistinguishable from that of the EB scheme for both $L^2$ and $H^1$ broken norms. 
Figure \ref{fig:Pois2D_pointError}a shows the absolute error of the solution obtained with the RB scheme for the case $\vert \mathcal{N}_{\Gamma_1} \vert=13$ and $\vert \mathcal{N}_{\Gamma_2} \vert=5$.
The same problem is solved considering also a curved interface (Figure \ref{fig:Pois2D_pointError}b), which implies the presence of geometrically non-conforming surfaces. Despite the gaps at the subdomain interface, the RB scheme is able to accurately transfer the solution across the interface.  
Error peaks can be noticed on the extremes nodes of the flat interface, due to the larger solution gradient along the horizontal direction. 

\begin{figure}
    \centering
    \includegraphics[width=\linewidth]{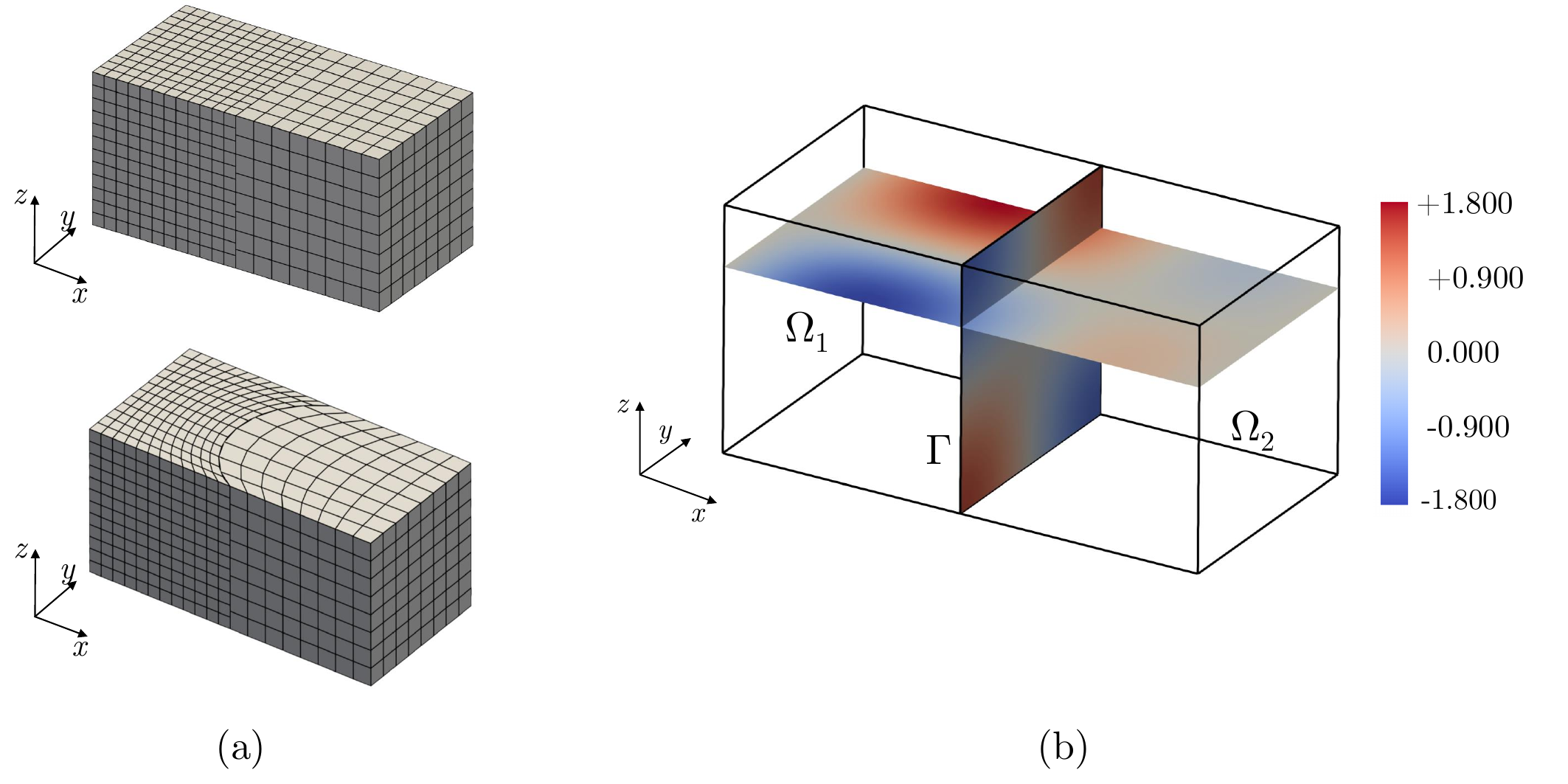}
    \caption{3D Poisson problem: (a) discretized domain with flat interface (top) and curved interface (bottom); (b) behavior of the manufactured solution.}
    \label{fig:Poisson3D_setting}
\end{figure}

\begin{figure}
    \centering
    \includegraphics[width=0.75\linewidth]{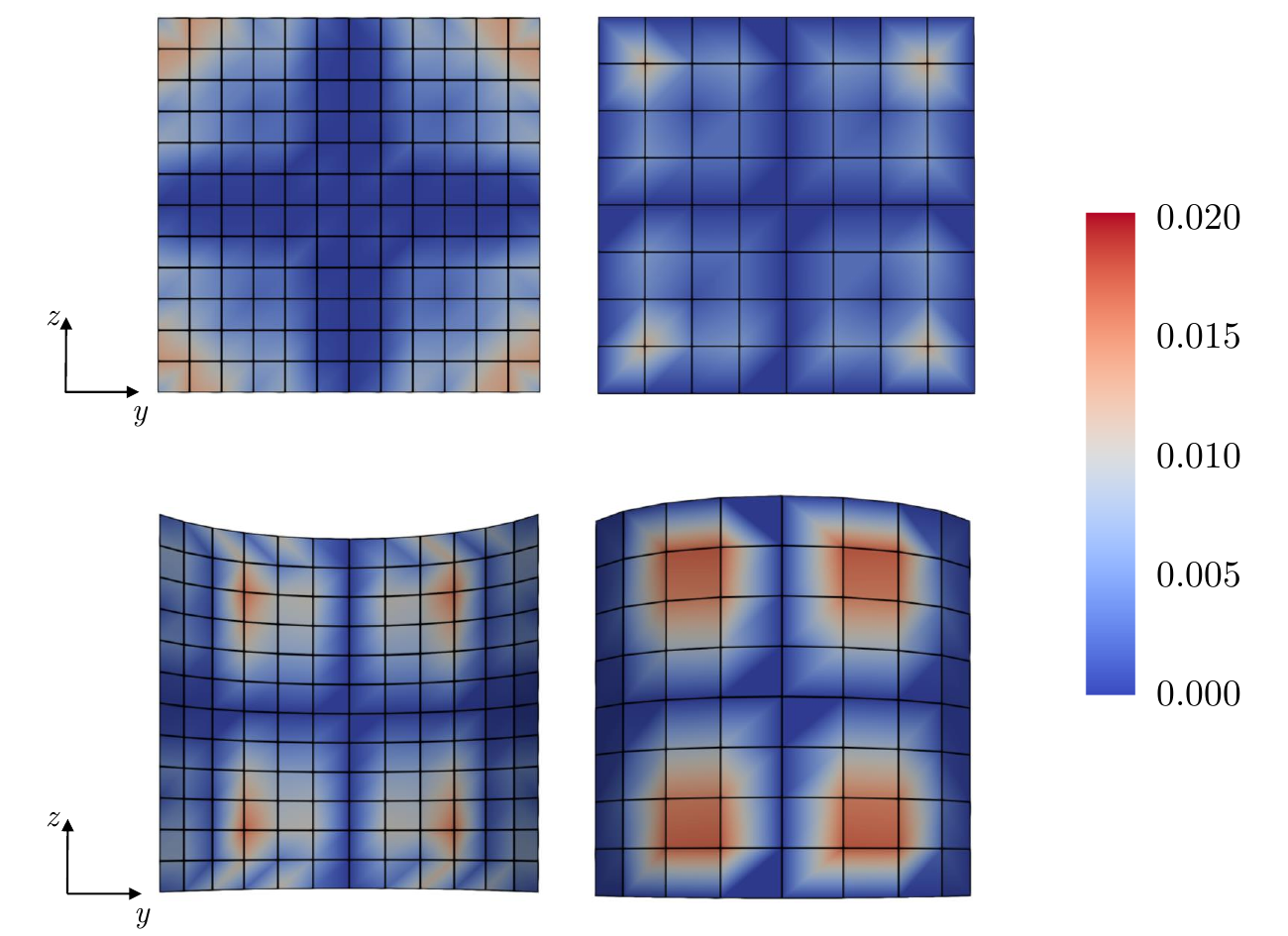}
    \caption{3D Poisson problem: absolute error of the manufactured solution on the interface for an example of discretization. Flat interface on top, curved interface on bottom.}
    \label{fig:Poisson3D_error}
\end{figure}

Consider now the 3D setting depicted in Figure \ref{fig:Poisson3D_setting}a, with $\Omega = [0,2] \times [0,1] \times [0,1]$. As usual, we partition $\Omega$ into two subdomains, $\Omega_1$ and $\Omega_2$, with either a flat or a curved interface $\Gamma$. Both subdomains are discretized by linear hexahedral elements with different characteristic size such that $h_1/h_2 = 2/3$.
The selected manufactured solution is $u(x,y,z) = \cos(\pi y) \cos(\pi z)(2x - x^2 + \sin(\pi x))$ (Figure \ref{fig:Poisson3D_setting}b).
The outcome of the
convergence study in $L^2$ and $H^1$ broken norms 
yields the same observations as for the 2D case (Figure \ref{Poisson3D_conv}).
Looking at the point-wise absolute errors obtained on $\Gamma$
(Figure \ref{fig:Poisson3D_error}), we can observe that, despite the slightly larger errors of the curved case, the accuracy proves to be satisfactory. Note also that the error distribution preserves the symmetry of the problem at hand.

%


\subsection{Realistic application: fluid-structure interaction problem} \label{sect5:Application}
\begin{figure}
    \centering
    \includegraphics[width=0.7\linewidth]{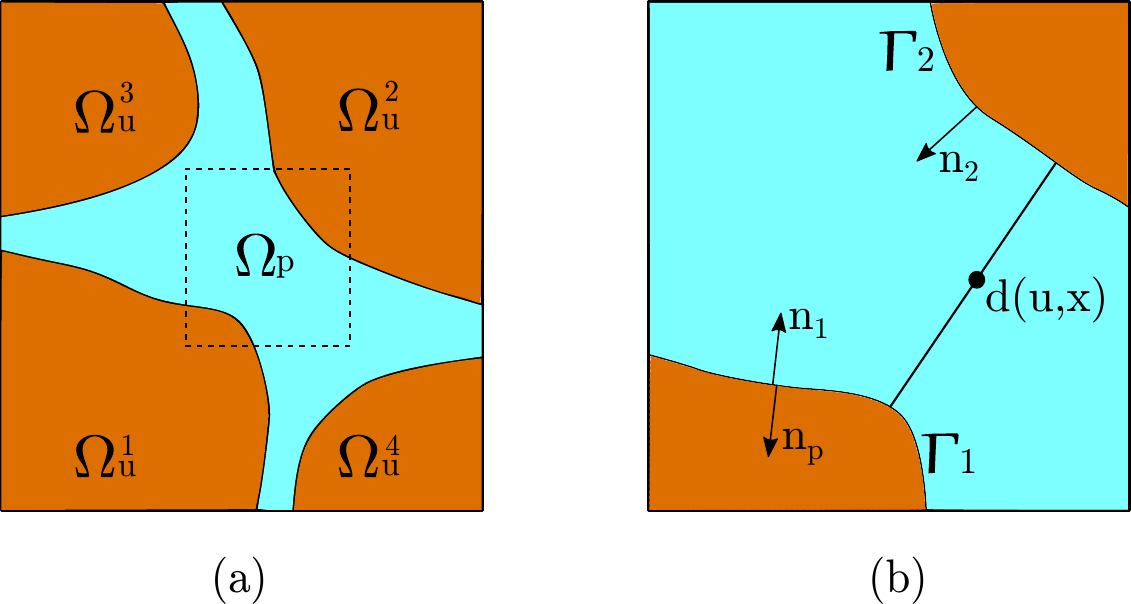}
    \caption{Sketch of FSI problem: (a) overall domain decomposed into $\Omega_p$ and $\Omega_u$; (b) zoom on the interface between two solid volumes and definition of the channel aperture $d$.}
    \label{fig:Pore_scale_sketch}
\end{figure}

The proposed mortar algorithm is finally used in a realistic test case simulating a fluid-structure interaction (FSI) problem. 
Let us consider the setting sketched in Figure \ref{fig:Pore_scale_sketch}, where the overall domain $\Omega \subset \mathbb R^3$ is split into two physical regions: 
\begin{itemize}
    \item a solid region $\Omega_u$, resulting from the union of $G$ volumes:
    \begin{align}
        \Omega_u = \bigcup_{i=1}^G \Omega_{u}^i,
    \end{align}
    where $\Omega_{u}^i$ denotes the domain portion occupied by each individual solid constituent. The symbol $\boldsymbol n_i$ denotes the outer normal vector to $\partial\Omega_u^i$, with $\partial\Omega_{D,u}^i$ and $\partial\Omega_{N,u}^i$ the corresponding Dirichlet and Neumann portions; 
    \item a fluid region $\Omega_p$, occupying the void space $\Omega\backslash\Omega_u$, with $\partial\Omega_{D,p}$ and $\partial\Omega_{N,p}$ its Dirichlet and Neumann boundaries and $\boldsymbol n_p$ the outer normal vector to $\partial \Omega_p$.
\end{itemize}
The selected setting can be representative of a FSI problem describing the motion of fluid in a porous material, where $\Omega_u$ are the solid grains and $\Omega_p$ the fluid channels.

$\Gamma_i$ denotes the interface between each solid volume $\Omega_u^i$ and the fluid region $\Omega_p$, with $\Gamma\equiv\cup_{i=1}^G\Gamma_i$ the union of all interfaces.
We assume an elastic behavior with infinitesimal strains for the solid region, while the fluid is incompressible and its flow in the channel is governed by a laminar regime. Body forces 
are neglected. The strong form of the considered boundary value problem reads: Given $\bar{\boldsymbol{t}}: \partial \Omega_{N,u} \rightarrow \mathbb{R}^3$ and $\bar p: \partial \Omega_{D,p} \rightarrow \mathbb{R}$, find $\boldsymbol{u}: \Omega_u \rightarrow \mathbb{R}^3$ and $p: \Omega_p\rightarrow \mathbb{R}$ such that:

\begin{subequations}
\begin{align}
    -\nabla \cdot \boldsymbol{\sigma}(\boldsymbol{u}) &= 0  \ \ \ \  \text{in} \ \ \Omega_u\ \  \ &\text{(linear momentum balance)} \\
    \nabla \cdot \boldsymbol q(\boldsymbol u,p) &= 0 \ \ \ \ \text{in} \ \ \Omega_p \ \ \  &\text{(mass balance)} \\
    \boldsymbol{\sigma}(\boldsymbol{u}) \cdot \boldsymbol n_i - p\boldsymbol n_p &= 0 \ \ \ \ \text{in} \ \ \Gamma \ \  \ &\text{(traction balance)} \label{eq: grain_pressure}\\
    \boldsymbol q \cdot \boldsymbol n_i &= 0 \ \ \ \ \text{in} \ \ \Gamma \ \ \ &\text{(interface flux)}\\
    \boldsymbol{u} \cdot \boldsymbol n_i &=  0 \ \ \ \, \text{in} \ \ \partial\Omega^i_{D,u} \ \ \ &\text{(boundary displacement)} \\
    \boldsymbol{\sigma}(\boldsymbol{u}) \cdot \boldsymbol n_i &= \bar{\boldsymbol{t}} \ \ \ \ \text{in} \ \ \partial \Omega_{N,u} \ \ \ &\text{(boundary traction)} \\
    p &= \bar p \ \ \ \ \text{in} \ \ \partial\Omega_{D,p} \ \ \ &\text{(boundary pressure)} \\
    \boldsymbol q \cdot \boldsymbol n_p &= 0 \ \ \ \ \text{in} \ \ \partial\Omega_{N,p} \ \ \ &\text{(boundary flux)}
\end{align}
\label{eq:IBVP_porescale}
\end{subequations}   
%
The following constitutive relationships close the problem:
\begin{itemize}
    \item $\boldsymbol{\sigma}(\boldsymbol{u}) = \boldsymbol{\mathsf{C}}:\nabla^s\boldsymbol u$ is the Cauchy stress tensor, with $\boldsymbol{\mathsf{C}}$ the fourth-order elasticity tensor and $\nabla^s$ the symmetric gradient operator;
    \item $\boldsymbol q(\boldsymbol u,p) = -C_f(\boldsymbol u)\mu^{-1}\nabla p$ is the volumetric flux in the fluid domain, varying linearly with the pressure gradient according to the laminar flow hypothesis \cite{witherspoon1980validity}. The conductivity $C_f$ is defined as \cite{garipov2016discrete}:
    \begin{align}
        C_f = \frac{d(\boldsymbol u)^3}{12},
        \label{eq:conductivity_porescale}
    \end{align}
    with $d(\boldsymbol u)$ representing the channel aperture of the void space (see Figure \ref{fig:Pore_scale_sketch}b).   
\end{itemize}
\begin{figure}
    \centering
    \includegraphics[width=0.7\linewidth]{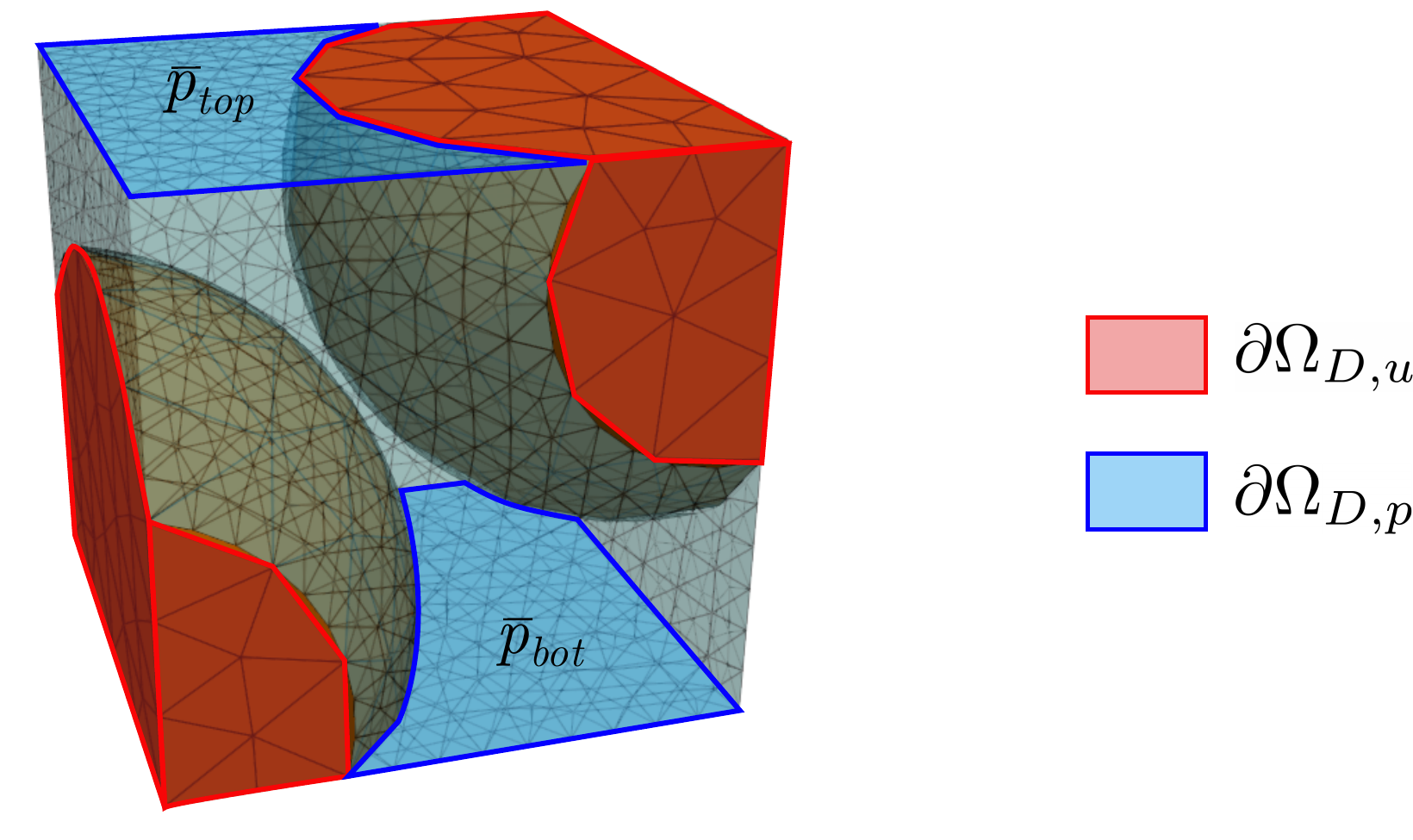}
    \caption{FSI problem: non-conforming grid.}
    \label{fig:pore_scale_grid}
\end{figure}
\begin{table}
    \centering
    \begin{tabular}{l l l}
    \toprule
    Quantity & Value & Unit \\[2pt]
    \midrule
    Young's modulus ($E$) & 1 $\times$ $\mathrm{10^{5}}$ & [kPa] \\[2pt]
    Poisson's ratio ($\nu$) & 0.3 & [-]\\[2pt]
    Fluid viscosity ($\mu$) & 1 $\times$ $\mathrm{10^{-6}}$ & [kPa$\cdot$s]\\[2pt]
    Domain size $x-y-z$ & 500 & $\mathrm{\mu m}$\\[2pt]
    Pressure at top $(\bar p_{top}$) & 5.1 $\times$ $\mathrm{10^{3}}$ & [kPa]\\[2pt]
    Pressure at bottom ($\bar p_{bot}$) & 5.0 $\times$ $\mathrm{10^{3}}$ & [kPa]\\[2pt]
    \hline
    \end{tabular}
    \caption{Material properties and domain dimensions of the FSI problem.}
    \label{tab:pore_scale}
\end{table}
Problem \eqref{eq:IBVP_porescale} is solved taking $G=2$ in the box-shaped domain whose discretization is shown in Figure \ref{fig:pore_scale_grid}. Material properties and domain dimensions are summarized in Table \ref{tab:pore_scale}. A pressure gradient is imposed by setting the pressures at top and bottom of the fluid domain equal to $\bar p_{top}$ and $\bar p_{bot}$, respectively. 
No body traction is considered.

The strong problem \eqref{eq:IBVP_porescale} is solved by a standard finite element method with linear elements, where different and independent discretizations are introduced for $\Omega_u$ and $\Omega_p$. The use of non-conforming grids for the solid and fluid region of the domain allows for a simpler and more flexible process of mesh generation and refinement.
A sequential algorithm is employed to address the coupling between the fluid pressure in $\Omega_p$ and the solid displacement in $\Omega_u$. The mass balance equation is solved first, then fluid pressure is transferred on the solid side by the RB mortar algorithm and the linear momentum balance is solved to compute the solid displacement field. 
The conductivity $C_f$ is updated by equation \eqref{eq:conductivity_porescale}, where $d(\boldsymbol u)$ is obtained by transferring $\boldsymbol u$ from the solid to the fluid side of $\Gamma$ again with the RB mortar algorithm. Notice that the master and the slave side are swapped when transferring $p$ and $\boldsymbol{u}$. 

\begin{figure}
    \centering
    \includegraphics[width=0.6\linewidth]{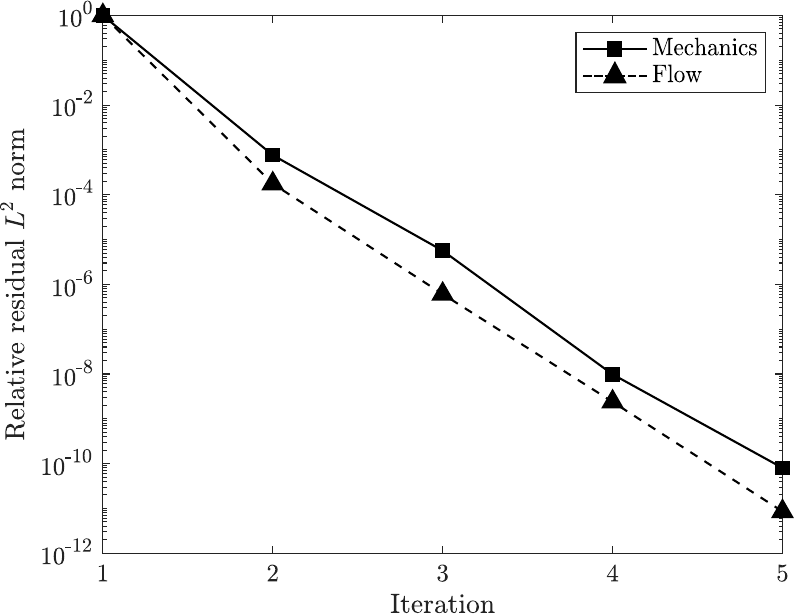}
    \caption{FSI problem: convergence profiles for flow ($p$) and mechanical ($\boldsymbol{u}$) problems.}
    \label{fig:conv_poreScale}
\end{figure}

\begin{figure}
\centering
\subfloat[]{\includegraphics[width = 0.51 \textwidth]{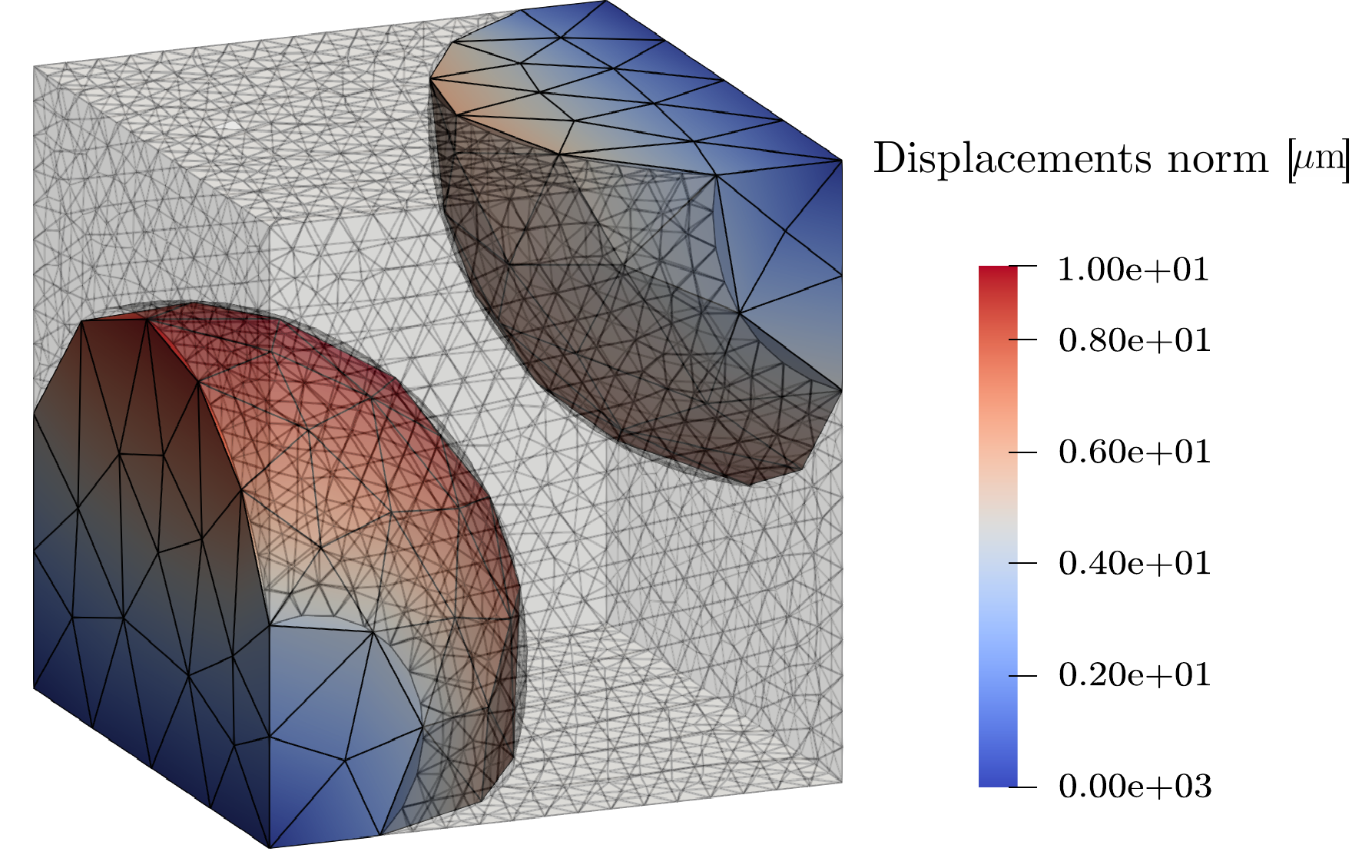}\label{fig:res_grains}}
\hspace{0.25cm}
\subfloat[]{\includegraphics[width = 0.45 \textwidth]{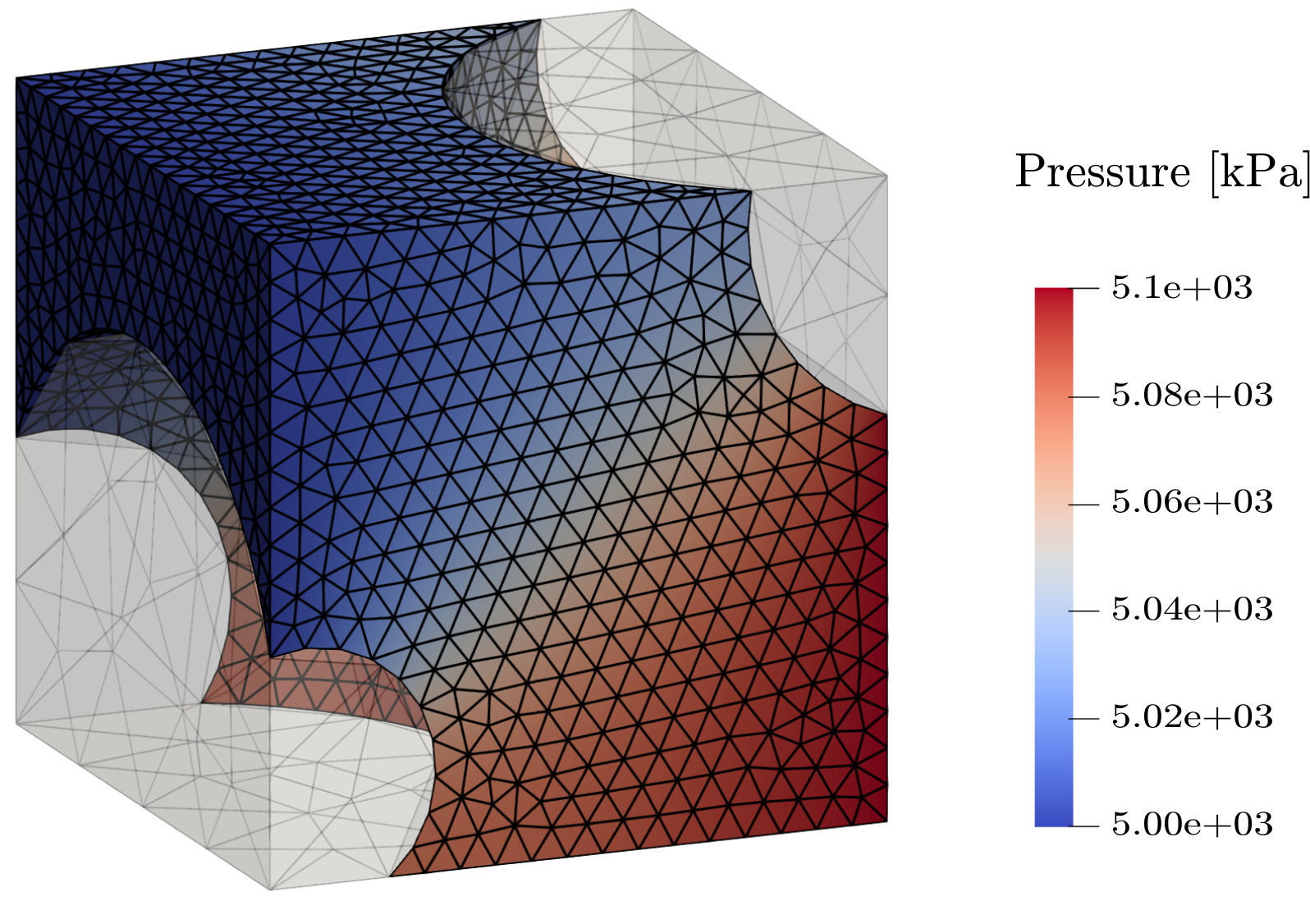}\label{fig:res_fluid}} 
\caption{FSI problem: contour plot of the converged solution. (a) Solid displacement norm; (b) fluid pressure field.}
\label{Fig:Pore_scale_results}
\end{figure}

The entire procedure is repeated iteratively until convergence, which is achieved when the following conditions are met:
\begin{align*}
    \frac{\|\boldsymbol{u}^{k+1}-\boldsymbol u^k \|}{\|\boldsymbol u^k\|} \leq \epsilon, \ \ \ \ \ \ \ \ \  \frac{\| p^{k+1}-p^k \|}{\|p^k\|} \leq \epsilon,
\end{align*}
with $k$ the iteration counter and $\epsilon$ a given tolerance, set to $10^{-9}$ for this case.
The convergence profiles for the $L^2$ norms of $\boldsymbol{u}$ and $p$ are reported in Figure \ref{fig:conv_poreScale}, showing a linear behavior as expected. 
Figure \ref{Fig:Pore_scale_results} provides the resulting contours of the pressure field and the displacement norm.

\section{Conclusions}



In this work, we have introduced a novel algorithm for the efficient computation of the interface integral appearing when using the mortar method 
in non-conforming discretizations. 
%
The proposed algorithm leverages the RBF-based mesh-free interpolation to replace the projection operator needed to transfer the basis functions defined on the master side of a non-conforming interface onto the slave side.
The new procedure provides
a more efficient alternative to standard geometric projections, replacing the solution to a set of non-linear problems with a simple interpolated function evaluation. The proposed mesh-free approach significantly streamlines implementation issues, especially across general three-dimensional settings. Moreover, we have employed a rescaling procedure for the interpolant that ensures that the necessary conditions for the consistency of the interface constraint are met. We also propose an appropriate combination of interpolation parameters allowing for an effective outcome of the mortar algorithm while minimizing the potential burden for the final user.

A set of numerical tests have been carried out in both 2D and 3D to analyze the method accuracy and convergence, demonstrating that the basic properties of the standard mortar approach are preserved. 
Results from geometrically non-conforming test cases and fluid-structure interaction problems show the applicability of the proposed approach in realistic settings and validate its flexibility and robustness.

\vspace{0.5cm}
\noindent {\bf Acknowledgements.} M. Ferronato was supported by the RESTORE (REconstruct subsurface heterogeneities and quantify Sediment needs TO improve the REsilience of Venice saltmarshes) PRIN 2022 PNRR Project, Funded by the European Union - Next Generation EU, Mission 4, Component 1 CUP MASTER B53D2303363001.









\bibliography{mybibfile}

\end{document}